\documentclass[11pt]{amsart}

\usepackage{amsfonts,amsmath,amssymb,amsthm}
\allowdisplaybreaks[4]

\usepackage[colorlinks,linkcolor=red,citecolor=green]{hyperref}
\usepackage{makecell}
\usepackage{graphics}
\usepackage{tikz}
\newcommand{\circled}[1]{
  \tikz[baseline=(char.base)]{
    \node[shape=circle,draw,inner sep=1pt] (char) {#1};
  }
}

\newtheorem{theorem}{Theorem}
\newtheorem{lemma}[theorem]{Lemma}
\newtheorem{remark}[theorem]{Remark}

\newtheorem{corollary}[theorem]{Corollary}

\numberwithin{theorem}{section}
\numberwithin{equation}{section}

\begin{document}

\setlength{\baselineskip}{1.1\baselineskip}

\title[The log-concavity of eigenfunction to complex Monge-Amp\`{e}re operator]{The log-concavity of eigenfunction to complex Monge-Amp\`{e}re operator in $\mathbb{C}^2$}

\author{Wei Zhang}
\address{School of Mathematics and Statistics\\
Lanzhou University\\
Lanzhou, 730000, Gansu Province, China.}
\email{zhangw@lzu.edu.cn}

\author{Qi Zhou}
\address{School of Mathematics and Statistics\\
Lanzhou University\\
Lanzhou, 730000, Gansu Province, China.}
\email{zhouqi2025@lzu.edu.cn, zhouqimath20@lzu.edu.cn}

\maketitle

\begin{abstract}
 Following the authors' recent work \cite{Zhang-Zhou2025}, we further explore the convexity properties of solutions to the Dirichlet problem for the complex Monge-Amp\`{e}re operator. In this paper, we establish the $\log$-concavity of solutions to the Dirichlet eigenvalue problem for the complex Monge-Amp\`{e}re operator on bounded, smooth, strictly convex domain in $\mathbb{C}^2$. Our approach combines a constant rank theorem adapted to the study of real convexity in the complex setting with a deformation method.
\end{abstract}



\section{Introduction}

\subsection{Main result}
The classical Dirichlet eigenvalue problem
\[\begin{cases}
 \Delta u=-\lambda u & \mathrm{in} \ \Omega,\\
 \hspace{3mm}u=0 & \mathrm{on} \ \partial\Omega
\end{cases}\]
has been widely studied, where $\Omega\subset\mathbb{R}^n$ is a bounded convex domain. Among these studies, the investigation of the $\log$-concavity of its solutions, also known as the  eigenfunctions of the Laplacian operator, is of particular importance. This property has significant applications in geometry, physics, and other fields. Initially, in 1976, Brascamp and Lieb \cite{Brascamp-Lieb1976} employed parabolic methods to prove the $\log$-concavity of the first Dirichlet eigenfuction of the Laplacian operator of a bounded convex domain, leading to the Brunn-Minkowski inequality for the first eigenvalue as following:
\[\lambda\left((1-t)\Omega_0+t\Omega_1\right)^{-1/2}\geq(1-t)\lambda(\Omega_0)^{-1/2}+t\lambda(\Omega_1)^{-1/2},\]
where $t\in[0,1]$, $\Omega_0, \Omega_1$ are nonempty convex bodies in $\mathbb{R}^n$. In 1985, Singer, Wong, Yau and Yau \cite{Singer-Wong-Yau-Yau1985} utilized the maximum principle combined with the continuity method to obtain the same concavity result and provided estimates for the gap between the first and second non-zero eigenvalues. The optimal lower bound, known as the famous fundamental gap conjecture, was completely resolved by Andrews and Clutterbuck \cite{Andrews-Clutterbuck2011}.

The eigenvalue problem for elliptic operators plays a significant role in the study of partial differential equations (PDEs). In addition to the previously mentioned eigenvalue problem for the Laplacian operator, the eigenvalue problem for fully nonlinear elliptic operators has also been extensively studied by many authors. Lions \cite{Lions1985} first investigated the eigenvalue problem for the Monge-Amp\`{e}re operator on uniformly convex domain $\Omega$ with smooth boundary in $\mathbb{R}^n$. He proved the existence of a unique positive constant $\lambda(\Omega)$ and a unique non-zero convex solution $u\in C^\infty(\Omega)\cap C^{1,1}(\bar{\Omega})$ solving
\[\begin{cases}
 \det (u_{ij})=\lambda(\Omega)(-u)^n & \mathrm{in}\ \Omega,\\
 \hspace{11mm}u=0 & \mathrm{on} \ \partial\Omega.
\end{cases}\]
Later, Wang \cite{Wang1994} generalized this result to the $k$-Hessian operator. More details on the eigenvalue problem for the Monge-Amp\`{e}re operator can be found in Chapter 11 of Le's recent book \cite{Le2024}.

In the complex setting, Kutev and Ramadanov \cite{Kutev-Ramadanov1989,Kutev-Ramadanov1990} were the first to initiate such studies. In 2023, Badiane and Zeriahi \cite{Badiane-Zeriahi2023} studied the existence and uniqueness of the first eigenvalue and an associated eigenfunction with Dirichlet conditions for the complex Monge-Amp\`{e}re operator on a bounded strongly pseudoconvex domain in $\mathbb{C}^n$. Soon after, Chu, Liu and McCleerey \cite{Chu-Liu-McCleerey-arXiv2024} extended this result to the complex $k$-Hessian operator on manifolds, proving the eigenfunction $u\in C^\infty(\Omega)\cap C^{1,1}(\bar{\Omega})$. In more detail, assume that $\Omega\subset\mathbb{C}^n$ is a bounded strongly pseudoconvex domain with smooth boundary, and that $0<\psi\in C^\infty(\bar{\Omega})$. Then there exists a real number $\lambda(\Omega,\psi)$ and a plurisubharmonic (PSH) function $u\in C^\infty(\Omega)\cap C^{1,1}(\bar{\Omega})$ such that the pair $(\lambda,u)$ is a unique solution to the eigenvalue problem
\begin{equation}\label{Equation-Dirichletproblem}
\begin{cases}
\det(u_{i\bar{j}})=\lambda(\Omega)\psi(-u)^n & \mathrm{in}\ \Omega,\\
\hspace{11mm}u=0 & \mathrm{on} \ \partial\Omega.
\end{cases}
\end{equation}

The aim of this paper is to study the $\log$-concavity of eigenfunction for the complex Monge-Amp\`{e}re operator, specifically to prove that the solution to the Dirichlet probelm \eqref{Equation-Dirichletproblem} is $\log$-concave when $\psi=1$ and $n=2$. We state the theorem as follows.
\begin{theorem}\label{Theorem-Main}
Let $\Omega$ be a bounded, smooth and strictly convex domain in $\mathbb{C}^2$. Suppose that $u\in\mathrm{PSH}(\Omega)\cap C^\infty(\Omega)\cap C^{1,1}(\bar{\Omega})$ is the unique solution of Dirichlet problem
\begin{equation}\label{Equ:EVP-u}
\begin{cases}
\det(u_{i\bar{j}})=\lambda(\Omega)(-u)^2 & \mathrm{in}\ \Omega,\\
\hspace{11mm}u<0 & \mathrm{in}\ \Omega,\\
\hspace{11mm}u=0 & \mathrm{on}\ \partial\Omega.
\end{cases}
\end{equation}
Then the function $v=-\log(-u/4)$ is strictly convex in $\Omega$.
\end{theorem}
\begin{remark}
An earlier version of this work was posted on arXiv as 	arXiv:2505.12817 and was subsequently withdrawn. During the revision of this manuscript, Chen, Li and Ma \cite{Chen-Li-Ma2026} obtained a more general result. More precisely, they established the strict real $\log$-concavity of the first eigenfunction of the complex $2$-Hessian operator in arbitrary complex dimension and further derived a Brunn-Minkowski inequality for its first eigenvalue. In complex dimension $2$, their $\log$-concavity result coincides with Theorem \ref{Theorem-Main}.
\end{remark}

\subsection{Previous works}
Convexity has been a central topic in PDEs for a long time. In this subsection, we briefly review the historical progress on the convexity of solutions to elliptic PDEs, emphasizing the $\log$-concavity of eigenfunctions for elliptic operators. For more detailed discussions, we refer the reader to the authors' recent article \cite{Zhang-Zhou2025}. 

In 1971, Makar-Limanov \cite{Makar-Limanov1971} studied the Dirichlet problem for the Laplace equation
\[\begin{cases}
 \Delta u=-1 & \mathrm{in}\ \Omega,\\
 \hspace{3mm}u=0 & \mathrm{on}\ \partial\Omega,
\end{cases}\]
where $\Omega\subset\mathbb{R}^2$ is a bounded, strictly convex domain with smooth boundary, and established the strict concavity of the function $\sqrt{u}$ on $\Omega$. His main idea was to introduce the auxiliary function
\[P=v^4\det(\nabla^2v),\]
where $v=\sqrt{u}$ and $u$ is a solution of the Dirichlet problem, and then demonstrate that $P$ is superharmonic in $\Omega$. This approach is known as the $P$-function method. Later, Acker, Payne and Philippin \cite{Acker-Payne-Philippin1981} applied this method to provide an alternative proof of the result by Brascamp and Lieb, namely, to establish the $\log$-concavity of the first Dirichlet eigenfunction for the Laplacian operator on planar convex domains. More recently, Ma, Shi and Ye \cite{Ma-Shi-Ye2012}, Jia, Ma and Shi \cite{Jia-Ma-Shi2023}, advanced this technique to derive similar results in higher dimensions.

A breakthrough in the study of convexity properties of solutions to linear and quasilinear elliptic equations was achieved by Korevaar \cite{Korevaar-1983,Korevaar1983} in the early 1980s. For a bounded convex domain $\Omega\subset\mathbb{R}^n$ and a solution $u$ , he introduced the convexity test function
\[\mathcal{C}(x,y,\mu)=u\left((1-\mu)x+\mu y\right)-(1-\mu)u(x)-\mu u(y)\]
for $(x,y,\mu)\in\bar{\Omega}\times\bar{\Omega}\times[0,1]$, where the nonnegativity of $\mathcal{C}(x,y,\mu)$ implies the concavity of $u$. This method, now widely known as Korevaar's concavity maximum principle, has become a fundamental tool in the field; see \cite{Kennington1985,Kawohl1986,Sakaguchi1987,Langford-Scheuer2021,Borrelli-Mosconi-Squassina2024,Mosconi-Riey-Squassina2024}. Among its applications, Korevaar \cite{Korevaar1983} used this method to rederive the $\log$-concavity of the first Dirichlet eigenfunction for the Laplacian operator, as also discussed in Caffarelli and Spruck \cite{Caffarelli-Spruck1982}.

In 1997, Alvarez, Lasry and Lions \cite{Alvarez-Lasry-Lions1997} developed a new approach, the convex envelope method, to establish the convexity of viscosity solution of general fully nonlinear degenerate elliptic equation
\[F(\nabla^2v,\nabla v,v,x)=0\]
in a bounded convex domain $\Omega$ of $\mathbb{R}^n$. In particular, they provided a new proof of the concavity result originally established by Brascamp and Lieb \cite{Brascamp-Lieb1976}. Naturally, this also serves as a powerful tool for studying convexity problems, as discussed in \cite{Bianchini-Salani2013,Crasta-Fragala2015,Ishige-Nakagawa-Salani2016}.

Another significant method for studying the convexity of solutions to PDEs is the constant rank theorem. This approach was initially introduced by Caffarelli and Friedman \cite{Caffarelli-Friedman1985} for semilinear elliptic equations in two dimensions. By combining it with the deformation process, they were able to reprove Makar-Limanov's result as well as the $\log$-concavity of the first Dirichlet eigenfunction for the Laplacian operator in the two dimensional case. In higher dimensions, Korevaar and Lewis \cite{Korevaar-Lewis1987} utilized the same approach to derive corresponding proofs.

As we mentioned earlier, in 1994, Wang \cite{Wang1994} studied the existence and uniqueness of the eigenvalue problem for the $k$-Hessian operator, specifically considering the Dirichlet problem
\begin{equation}\label{Equation-Dirichlet problem-sigma_k}
\begin{cases}
 \sigma_k(\nabla^2u)=\lambda(\Omega)(-u)^k & \mathrm{in} \ \Omega,\\
 \hspace{12mm}u=0 & \mathrm{on} \ \partial\Omega,
\end{cases}
\end{equation}
where $1<k<n$ and $\Omega\subset\mathbb{R}^n$ is a bounded domain satisfying certain convexity conditions. Later, Liu, Ma and Xu \cite{Liu-Ma-Xu2010} applied the constant rank theorem to prove the $\log$-concavity of solutions to the Dirichlet problem \eqref{Equation-Dirichlet problem-sigma_k} for the case $k=2$ and $n=3$, that is, the $\log$-concavity of eigenfunctions of the $2$-Hessian operator. They also established the Brunn-Minkowski inequality for $2$-Hessian eigenvalues. Subsequently, Salani \cite{Salani2012} provided a new proof of this result. More recently, Li, Ma and Salani \cite{Li-Ma-Salani2026} extend the result to arbitrary dimensions. In addition, several lower dimensional convexity results of solutions to fully nonlinear elliptic boundary value problems appeared via constant rank theorems, for example \cite{Huang2019,Zhang-Zhou2023,Chen-Jia-Xiong2025}.

Recently, we have discovered that the constant rank theorem is also effective for studying convexity problems in the complex setting. In \cite{Zhang-Zhou2025}, we considered the Dirichlet problem for the complex Monge-Amp\`{e}re equation
\[\begin{cases}
 \det(u_{i\bar{j}})=1 & \mathrm{in}\ \Omega,\\
 \hspace{1.1cm}u=0 & \mathrm{on}\ \partial\Omega,
\end{cases}\]
where $\Omega\subset\mathbb{C}^2$ is a bounded, smooth, and strictly convex domain, and proved that $-\sqrt{-u/2}$ is strictly convex. More recently, Chen, Hu and Sheng \cite{Chen-Hu-Sheng2026} obtained the same result by a different method. 

\subsection{Remarks on the proof}\label{Section 1.3}
The proof of Theorem \ref{Theorem-Main} is based on the same general framework as our previous work \cite{Zhang-Zhou2025}, which is inspired by the classical methods of Caffarelli and Friedman \cite{Caffarelli-Friedman1985} and Korevaar and Lewis \cite{Korevaar-Lewis1987}. The proof consists of two main steps: deriving convexity estimates near the boundary and establishing a constant rank theorem for the transformed equation satisfied by $v=-\log(-u/4)$, where $u$ solves the complex Monge-Amp\`{e}re Dirichlet problem \eqref{Equ:EVP-u}.

Due to the strict convexity of the domain $\Omega$, we can easily obtain the convexity estimates for $v$ near the boundary $\partial\Omega$. For details, see Lemma \ref{Lemma-BCE}.

The principal difficulty of this paper is the establishment of the constant rank theorem. The general strategy is to choose a suitable auxiliary function, derive an appropriate differential inequality for it, and then apply the strong maximum principle. There are two main difficulties in carrying out this argument. The first concerns the choice of the auxiliary function. In proving a constant rank theorem for solutions of general fully nonlinear elliptic equations $$F(\nabla^2v,\nabla v,v,x)=0\quad\mathrm{in}\ \Omega\subset\mathbb{R}^n,$$ Bian and Guan \cite{Bian-Guan2009} introduced the auxiliary function $$\phi(x)=\sigma_{l+1}(\nabla^2v(x))+q(\nabla^2v(x)),$$ where $$q(\nabla^2v(x))=\begin{cases}
\frac{\sigma_{l+2}(\nabla^2v(x))}{\sigma_{l+1}(\nabla^2v(x))} & \mathrm{if}\ \sigma_{l+1}(\nabla^2v(x))>0,\\
0 & \mathrm{if}\ \sigma_{l+1}(\nabla^2v(x))=0.
\end{cases}$$ Here, $l$ denotes the minimum rank of $\nabla^2v$ in $\Omega$. See Section \ref{Section 2} for more discussions on the elementary symmetric functions. In the present complex setting, however, this auxiliary function is not well suited to our problem in its original form. Indeed, unitary transformations naturally diagonalize Hermitian matrices, whereas the real convexity considered here is governed by the real Hessian $\nabla^2v$, which is generally not diagonalized by such transformations. Consequently, this makes our computations considerably more complicated. We therefore modify the construction of Bian and Guan and introduce auxiliary functions adapted to the complex structure; see Section \ref{Section 2}. Further motivation and details can be found in Section 1.3 and 2.2 of \cite{Zhang-Zhou2025}. The second major difficulty is to establish the required differential inequality for the modified auxiliary functions. This is the central task of this paper and involves highly nontrivial computations.

The rest of this paper is organized as follows. In Section \ref{Section 2}, we briefly introduce the preliminaries needed for the subsequent analysis, including the real and complex coordinate settings, the notation used throughout the paper, the choice of suitable coordinates, the construction of the auxiliary functions, and the corresponding computations. In Section \ref{Section 3}, we consider the two possible cases separately and establish the required differential inequality in each case, thereby proving the constant rank theorem. This section constitutes the main part of the paper. Finally, in Section \ref{Section 4}, we prove Theorem \ref{Theorem-Main}.

\section{Notation and Preliminaries}\label{Section 2}

In this section, we introduce the notation and preliminary facts that will be used throughout the paper. We consider the complex Monge-Amp\`{e}re equation
\begin{equation}\label{Equation-complexMA(eigenvalue)}
\det(u_{i\bar{j}})=\lambda(\Omega)(-u)^2\quad \mathrm{in}\ \Omega\subset\mathbb{C}^2.
\end{equation}
To avoid confusion with the notation for real derivatives used later, we shall express the complex derivatives using the coordinates $z_i$ and $\bar{z}_i$. Under the logarithmic transformation $v=-\log(-u/4)$, namely $u=-4e^{-v}$, we have $$u_{z_i}=4e^{-v}v_{z_i},\quad u_{z_i\bar{z}_j}=4e^{-v}(v_{z_i\bar{z}_j}-v_{z_i}v_{\bar{z}_j}),$$ and therefore the equation \eqref{Equation-complexMA(eigenvalue)} is transformed into
\begin{equation}\label{Equ-complexMA-v}
\det(v_{z_i\bar{z}_j}-v_{z_i}v_{\bar{z}_j})=\lambda.
\end{equation}

Since we study the real convexity of solutions, the relevant object is the real Hessian $\nabla^2v$. We therefore rewrite the equation \eqref{Equ-complexMA-v} in real coordinates. Let $z_i=x_i+\sqrt{-1}y_i$. For a real-valued function $v$, the complex derivatives are given by $$v_{z_i}=\frac{\partial v}{\partial z_i}=\frac{1}{2}\left(\frac{\partial v}{\partial x_i}-\sqrt{-1}\frac{\partial v}{\partial y_i}\right)=\frac{1}{2}(v_{x_i}-\sqrt{-1}v_{y_i})$$ and $$v_{\bar{z}_i}=\frac{\partial v}{\partial\bar{z}_i}=\frac{1}{2}\left(\frac{\partial v}{\partial x_i}+\sqrt{-1}\frac{\partial v}{\partial y_i}\right)=\frac{1}{2}(v_{x_i}+\sqrt{-1}v_{y_i}).$$ Moreover
\begin{equation}\label{Equality-v_{ibarj}}
v_{z_i\bar{z}_j}=\frac{\partial^2v}{\partial z_i\partial\bar{z}_j}=\frac{1}{4}(v_{x_ix_j}+v_{y_iy_j})+\frac{\sqrt{-1}}{4}(v_{x_iy_j}-v_{y_ix_j})
\end{equation}
and $$v_{z_iz_j}=\frac{\partial^2v}{\partial z_i\partial z_j}=\frac{1}{4}(v_{x_ix_j}-v_{y_iy_j})-\frac{\sqrt{-1}}{4}(v_{x_iy_j}+v_{y_ix_j}).$$

Let $A=(v_{z_i\bar{z}_j})$, $B=(v_{z_iz_j})$, and write the real Hessian of $v$ as $$\nabla^2v=\begin{pmatrix}
U & V\\
V^T & W
\end{pmatrix},$$ where $U=(v_{x_ix_j})$, $V=(v_{x_iy_j})$, and $W=(v_{y_iy_j})$, with $V^T$ denoting the transpose of $V$. Then 
\begin{equation}\label{Matrix-A}
A=\frac{1}{4}(U+W)+\frac{\sqrt{-1}}{4}(V-V^T),
\end{equation} 
and 
\begin{equation}\label{Matrix-B}
B=\frac{1}{4}(U-W)-\frac{\sqrt{-1}}{4}(V+V^T).
\end{equation}
Since $u$ is strictly plurisubharmonic and $v<0$, the transformation $ u_{z_i\bar{z}_j}=4e^{-v}(v_{z_i\bar{z}_j}-v_{z_i}v_{\bar{z}_j})$ implies that $A$ is positive definite. In particular, $A$ is invertible.
 
For convenience, we adopt the following notation for the real derivatives
\begin{align*}
&v_{x_1}=v_1, \ v_{x_2}=v_2, \ v_{y_1}=v_3, \ v_{y_2}=v_4, \\
&v_{x_1x_1}=v_{11}, \ v_{y_1y_1}=v_{33}, \ v_{x_1y_1}=v_{13}, \ \ldots
\end{align*}
Substituting \eqref{Equality-v_{ibarj}} into \eqref{Equ-complexMA-v} yields the real equation
\begin{equation}\label{Equation-real equation(v)}
\begin{aligned}
\lambda
=&~(v_{11}+v_{33})(v_{22}+v_{44})-(v_{12}+v_{34})(v_{21}+v_{43})-(v_{14}-v_{32})(v_{41}-v_{23})\\
&-(v_1^2+v_3^2)(v_{22}+v_{44})-(v_2^2+v_4^2)(v_{11}+v_{33})\\
&+(v_1v_2+v_3v_4)(v_{12}+v_{21}+v_{34}+v_{43})\\
&+(v_1v_4-v_2v_3)(v_{14}-v_{32}+v_{41}-v_{23}). \end{aligned}
\end{equation}
Thus $v$ satisfies the Dirichlet problem
\begin{equation}\label{Equation-Dirichlet problem-real}
\begin{cases}
F(\nabla^2v,\nabla v)=0 & \mathrm{in}\ \Omega,\\
\hspace{18mm} v\rightarrow 0 & \mathrm{on}\ \partial\Omega,
\end{cases}
\end{equation}
where
\begin{align*}
F(\nabla^2v,\nabla v)
=&~(v_{11}+v_{33})(v_{22}+v_{44})-(v_{12}+v_{34})(v_{21}+v_{43})-(v_{14}-v_{32})(v_{41}-v_{23})\\
&-(v_1^2+v_3^2)(v_{22}+v_{44})-(v_2^2+v_4^2)(v_{11}+v_{33})\\
&+(v_1v_2+v_3v_4)(v_{12}+v_{21}+v_{34}+v_{43})\\
&+(v_1v_4-v_2v_3)(v_{14}-v_{32}+v_{41}-v_{23})-\lambda.
\end{align*}
For simplicity, we denote $$F^{ij}=\frac{\partial F}{\partial v_{ij}},\ F^{v_k}=\frac{\partial F}{\partial v_k},\ F^{ij,kl}=\frac{\partial^2F}{\partial v_{ij}\partial v_{kl}},\ F^{ij,v_k}=\frac{\partial^2F}{\partial v_{ij}\partial v_k},\ F^{v_k,v_l}=\frac{\partial^2F}{\partial v_k\partial v_l}.$$

We shall follow the framework developed in our previous work \cite{Zhang-Zhou2025} to establish the constant rank theorem in this present setting. In preparation for the proof, we present the selection of appropriate coordinate systems, the construction of auxiliary functions, and the necessary computations related to these functions that will be used in the subsequent analysis. For the detailed motivation, calculations, and explanations, we refer the reader to Section 1.3, 2, and 3.1 of \cite{Zhang-Zhou2025}.

Define the auxiliary matrix
\begin{equation}\label{Auxiliary matrix K}
K=I-B\overline{A^{-1}}\overline{B}A^{-1},
\end{equation}
where $I$ denotes the identity matrix of size $2\times 2$, $A=(v_{z_i\bar{z}_j})_{1\leq i,j\leq 2}$ and $B=(v_{z_iz_j})_{1\leq i,j\leq 2}$. As proved in \cite{Zhang-Zhou2025}, the real Hessian $\nabla^2v$ and the matrix $K$ satisfy the rank identity
\begin{equation}\label{rank identity}
\mathrm{rank}(\nabla^2v)=\mathrm{rank}(K)+2.
\end{equation}
Consequently, the study of the constant rank property of $\nabla^2v$ is reduced to that of $K$.

For completeness, we recall the adapted coordinates system constructed in \cite{Zhang-Zhou2025}. By the Autonne-Takagi factorization, for any fixed point $x_0\in\Omega$, there exists an invertible complex matrix $P$ such that $$PAP^H=I,\quad PBP^T=\Lambda,$$ where $P^H$ denotes the conjugate transpose of $P$ and $\Lambda$ is a nonnegative diagonal matrix. These relations, together with the explicit forms of $A$ and $B$ given in \eqref{Matrix-A} and \eqref{Matrix-B}, imply that, after the complex linear coordinate transformation induced by $P$, the real Hessian $\nabla^2v$ is diagonal at $x_0$. Under this transformation, the auxiliary matrix $K$ satisfies $$\widetilde{K}=PKP^{-1}.$$ Hence, the eigenvalues of $K$, as well as the quantities $\mathrm{tr}(K)$ and $\det(K)$, are preserved. Moreover, after normalizing the transformation by $|\det(P)|=1$, the equation \eqref{Equ-complexMA-v} remains invariant. Thus, all local computations may be performed in this adapted coordinate system.

Define $$l=\min_{x\in\Omega}\mathrm{rank}(\nabla^2v(x)).$$
At a point $x_0\in\Omega$ where the minimal rank $l$ is achieved, the adapted coordinate system allows us to assume that $\nabla^2v(x_0)$ is diagonal. If $l=0$ or $l=1$, then the right-hand side of the real equation \eqref{Equation-real equation(v)} is nonpositive at $x_0$, while the left-hand side $\lambda$ is positive, which leads to a contradiction. If $l=4$, then $\nabla^2v$ is of full rank in $\Omega$, and the constant rank property is immediate. Hence, it suffices to consider the two cases $l=2$ and $l=3$. By the rank identity \eqref{rank identity}, these two cases correspond respectively to

$\bullet$ $\min_{x\in\Omega}\mathrm{rank}(K(x))=0$,

$\bullet$ $\min_{x\in\Omega}\mathrm{rank}(K(x))=1$.

The construction of the auxiliary functions is inspired by the ideas of Bian and Guan \cite{Bian-Guan2009}. Since the present equation differs from the original setting considered therein, we make suitable modifications to their construction. The motivation and details of these modifications have been discussed in Section 1.3 of our previous work \cite{Zhang-Zhou2025}. Accordingly, we choose the auxiliary functions 
\begin{equation}\label{Auxiliary function-l=2}
\phi(x)=\mathrm{tr}(K(x))+q(K(x))
\end{equation}
and $$\phi(x)=\det(K(x))$$ for the above two cases, respectively, where $$q(K(x))=\begin{cases}
\frac{\det(K(x))}{\mathrm{tr}(K(x))} & \mathrm{if}\ \mathrm{tr}(K(x))>0,\\
0 & \mathrm{if}\ \mathrm{tr}(K(x))=0.
\end{cases}$$ These auxiliary functions will be used in the subsequent sections to derive the differential inequalities required for the constant rank theorem.

To perform the differential computations of these auxiliary functions, we first write them in terms of the entries of the real Hessian $\nabla^2v$. In particular, the determinant of $K$ can be related to the fourth elementary symmetric function of $\nabla^2v$. We therefore recall the notation of $k$-th elementary symmetric functions. For notational convenience, we directly consider the general real dimension $n$. For $1\leq k\leq n$ and $\lambda=(\lambda_1,\lambda_2,\cdots,\lambda_n)\in\mathbb{R}^n$, define $$\sigma_k(\lambda)=\sum_{1\leq i_1<i_2<\cdots<i_k\leq n}\lambda_{i_1}\lambda_{i_2}\cdots\lambda_{i_k}.$$ We adopt the convention $\sigma_0(\lambda)=1$ and $\sigma_{-1}(\lambda)=\sigma_{n+1}(\lambda)=0$. For an $n\times n$ symmetric matrix $W$, we define $\sigma_k(W)=\sigma_k(\lambda(W))$, where $\lambda(W)=(\lambda_1(W),\lambda_2(W),\cdots,\lambda_n(W))$ are the eigenvalues of $W$. We also denote by $\sigma_{k-1}(W|i)$ the $(k-1)$-th elementary symmetric function obtained by setting $\lambda_i=0$, and by $\sigma_{k-2}(W|ij)$ ($i\neq j$) the $(k-2)$-th elementary symmetric function obtained by setting $\lambda_i=\lambda_j=0$.

The following derivative formulas for elementary symmetric functions will be frequently used in the subsequent calculations.
\begin{lemma}\label{Lemma-sigma_k}
Suppose that $W=(W_{ij})_{n\times n}$ is diagonal. For $k=1,2,\cdots,n$, then we have $$\frac{\partial\sigma_k(W)}{\partial W_{ij}}=\begin{cases}
\sigma_{k-1}(W|i) & \mathrm{if}\ i=j,\\
0 & \mathrm{if}\ i\neq j,
\end{cases}$$ and $$\frac{\partial^2\sigma_k(W)}{\partial W_{ij}\partial W_{pq}}=\begin{cases}
\sigma_{k-2}(W|ip) & \mathrm{if}\ i=j, p=q, i\neq p,\\
-\sigma_{k-2}(W|ip) & \mathrm{if}\ i=q, j=p, i\neq p,\\
0 & \mathrm{otherwise}.
\end{cases}$$
\end{lemma}

We now focus on the case where $\nabla^2v$ has minimal rank $2$. In this case, the auxiliary function is given by \eqref{Auxiliary function-l=2}. Using the definition of $K$ in \eqref{Auxiliary matrix K}, together with the standard real-complex correspondence $$4\begin{pmatrix}
A & B\\[6pt]
\overline{B} & \overline{A}
\end{pmatrix}=\begin{pmatrix}
I & -\sqrt{-1}I\\[6pt]
I & \sqrt{-1}I
\end{pmatrix}\begin{pmatrix}
U & V\\[6pt]
V^T & W
\end{pmatrix}\begin{pmatrix}
I & I\\[6pt]
\sqrt{-1}I & -\sqrt{-1}I
\end{pmatrix},$$ we have 
\begin{equation}\label{Equality-det K and det D^2v}
\det(K)=(\det(A))^{-2}\det(\nabla^2v),
\end{equation}
where $\det(A)$ is computed from \eqref{Matrix-A} as $$\det(A)=\frac{1}{16}\big((v_{11}+v_{33})(v_{22}+v_{44})-(v_{12}+v_{34})^2-(v_{14}-v_{32})^2\big).$$ Moreover, $$K=(\det(A))^{-2}\big((\det(A))^2I-B\overline{A^\ast}\overline{B}A^\ast\big),$$ where $A^\ast$ denotes the adjugate matrix of $A$. Hence, $$\mathrm{tr}(K)=(\det(A))^{-2}\mathrm{tr}\big((\det(A))^2I-B\overline{A^\ast}\overline{B}A^\ast\big).$$ Let $\psi=\mathrm{tr}\big((\det(A))^2I-B\overline{A^\ast}\overline{B}A^\ast\big)$. More explicitly, $\psi$ is given by
\begin{equation}\label{Equality-psi}
\begin{aligned}
\psi
=&~\frac{2}{16^2}\big((v_{11}+v_{33})(v_{22}+v_{44})-(v_{12}+v_{34})^2-(v_{14}-v_{32})^2\big)^2\\
&-\frac{1}{16^2}\bigg[\big((v_{11}-v_{33})(v_{22}+v_{44})+(v_{23}^2+v_{34}^2-v_{12}^2-v_{14}^2)\big)^2\\
&\hspace{1.2cm}+\big((v_{22}-v_{44})(v_{11}+v_{33})+(v_{14}^2+v_{34}^2-v_{12}^2-v_{23}^2)\big)^2\\[6pt]
&\hspace{1.2cm}+4\big(v_{12}v_{23}+v_{14}v_{34}-v_{13}(v_{22}+v_{44})\big)^2+4\big(v_{12}v_{14}+v_{23}v_{34}-v_{24}(v_{11}+v_{33})\big)^2\\[6pt]
&\hspace{1.2cm}+8\big(v_{12}v_{33}-v_{11}v_{34}+v_{13}(v_{14}-v_{32})\big)\big(v_{21}v_{44}-v_{22}v_{43}+v_{24}(v_{23}-v_{41})\big)\\
&\hspace{1.2cm}+8\big(v_{13}(v_{12}+v_{34})-v_{11}v_{32}-v_{33}v_{14}\big)\big(v_{24}(v_{21}+v_{43})-v_{23}v_{44}-v_{22}v_{14}\big)\bigg].
\end{aligned}
\end{equation}
Then, the auxiliary function in \eqref{Auxiliary function-l=2} can be rewritten as
\begin{equation}\label{Equality-phi}
\phi(x)=(\det(A))^{-2}\psi(x)+\frac{\sigma_4(\nabla^2v(x))}{\psi(x)}.
\end{equation}
For convenience, we continue to denote
\begin{equation}\label{q(D^2v(x))-new}
q(\nabla^2v(x))=\frac{\sigma_4(\nabla^2v(x))}{\psi(x)}.
\end{equation}

We next set up the notation for the case of minimal rank $2$ in the constant rank argument, following the framework developed by Bian and Guan \cite{Bian-Guan2009}. Suppose that $\nabla^2v$ attains its minimal rank $2$ at $x_0\in\Omega$. Let $\mathcal{O}$ be a sufficiently small neighborhood of $x_0$. For $x\in\mathcal{O}$, denote the eigenvalues of $\nabla^2v(x)$ by $\lambda_1(x)\geq\lambda_2(x)\geq\lambda_3(x)\geq\lambda_4(x)\geq 0$. There exists a positive constant $c$, independent of $x$, such that $\lambda_1(x)\geq\lambda_2(x)\geq c$. We define the “good" and “bad" index set by $G=\{1,2\}$, $B=\{3,4\}$. If there is no confusion, we also denote the “good" and “bad" eigenvalues of $\nabla^2v(x)$ by $G$ and $B$, respectively. Note that for any $\delta>0$, we may choose $\mathcal{O}$ small enough such that $\lambda_i(x)<\delta$ for all $i\in B$ and $x\in\mathcal{O}$. For two functions $h(x)$ and $k(x)$ defined on $\mathcal{O}$, we use the notation $h(x)=O(k(x))$ if there exists a positive constant $C$ such that $|h(x)|\leq C|k(x)|$.

At any fixed point $x\in\mathcal{O}$, we work in the coordinate system constructed previously, under which $\nabla^2v$ is diagonal and satisfies  
\begin{equation}\label{Equality-v_{11}+v_{33}=v_{22}+v_{44}}
v_{11}+v_{33}=v_{22}+v_{44}.
\end{equation}
Substituting the diagonal form of $\nabla^2v$ into the expression of $\psi$ in \eqref{Equality-psi}, we obtain $$\psi=\frac{4}{16^2}(v_{22}+v_{44})^2v_{11}v_{33}+\frac{4}{16^2}(v_{11}+v_{33})^2v_{22}v_{44}.$$ Since the good eigenvalues are uniformly bounded from above and below, we have $c\leq v_{22}\leq v_{11}\leq C$ for some positive constants $c$ and $C$ independent of $x\in\mathcal{O}$. Consequently, 
\begin{equation}\label{Estimate-psi}
\psi\geq\frac{4}{16^2}(v_{11}^2v_{22}v_{44}+v_{11}v_{22}^2v_{33})\geq c(v_{33}+v_{44}).
\end{equation}
On the other hand, it is also clear from the expression of $\psi$ that $\psi\leq C(v_{33}+v_{44})$. Furthermore, in view of \eqref{Equality-phi} and the uniform bounds for $v_{11}$ and $v_{22}$, we conclude that $c(v_{33}+v_{44})\leq\phi\leq C(v_{33}+v_{44})$. Therefore,
\begin{equation}\label{Estimate-v_{33}+v_{44}=phi}
0\leq c\phi\leq v_{33}+v_{44}\leq C\phi,
\end{equation}
and hence
\begin{equation}\label{Estimate-v_{33},v_{44}}
 v_{33}=O(\phi),\quad v_{44}=O(\phi).
\end{equation}
Using \eqref{Equality-v_{11}+v_{33}=v_{22}+v_{44}}, we further have
\begin{equation}\label{Estimate-v_{22}}
v_{22}=v_{11}+O(\phi).
\end{equation}

The following differential estimates for $q(\nabla^2v)$ in \eqref{q(D^2v(x))-new} were established in our earlier work \cite{Zhang-Zhou2025}. We recall them here for completeness.
\begin{lemma}[Lemma 3.2 in \cite{Zhang-Zhou2025}]\label{Lemma-first derivative of q}
Let $q(\nabla^2v)$ be defined as in \eqref{q(D^2v(x))-new}. Assume that $\nabla^2v$ is diagonal. Then
$$\frac{\partial q}{\partial x_\alpha}=O\bigg(\phi+\sum_{p,q\in B}|\nabla v_{pq}|\bigg).$$
\end{lemma}
\begin{lemma}[Proposition 3.4 in \cite{Zhang-Zhou2025}]\label{Lemma-q_{ij}}
Let $q(\nabla^2v)$ be defined as in \eqref{q(D^2v(x))-new}. Assume that $\nabla^2v(x)$ is diagonal for $x\in\mathcal{O}$. Then, for any $\alpha,\beta\in\{1,2,3,4\}$, we have
\begin{align*}
\frac{\partial^2q}{\partial x_\alpha\partial x_\beta}
=&~4\cdot 16v_{11}^{-1}\frac{v_{44}^2}{\sigma_1^2(B)}v_{33\alpha\beta}-8\cdot 16v_{11}^{-2}\frac{v_{44}^2}{\sigma_1^2(B)}v_{13\alpha}v_{13\beta}-8\cdot 16v_{11}^{-2}\frac{v_{44}^2}{\sigma_1^2(B)}v_{23\alpha}v_{23\beta}\\
&+4\cdot16v_{11}^{-1}\frac{v_{33}^2}{\sigma_1^2(B)}v_{44\alpha\beta}-8\cdot 16v_{11}^{-2}\frac{v_{33}^2}{\sigma_1^2(B)}v_{14\alpha}v_{14\beta}-8\cdot16v_{11}^{-2}\frac{v_{33}^2}{\sigma_1^2(B)}v_{24\alpha}v_{24\beta}\\
&-4\cdot16v_{11}^{-1}\frac{1}{\sigma_1^3(B)}(V_{3\alpha}V_{3\beta}+V_{4\alpha}V_{4\beta})-8\cdot16v_{11}^{-1}\frac{1}{\sigma_1(B)}v_{34\alpha}v_{34\beta}\\
&+O\bigg(\phi+\sum_{p,q\in B}|\nabla v_{pq}|\bigg),
\end{align*}
where $$V_{3\alpha}=v_{33\alpha}\sigma_1(B)-v_{33}\bigg(\sum_{i\in B}v_{ii\alpha}\bigg),\quad\mathrm{and}\ V_{4\alpha}=v_{44\alpha}\sigma_1(B)-v_{44}\bigg(\sum_{i\in B}v_{ii\alpha}\bigg).$$
\end{lemma}

We finish this section by recalling two lemmas of Bian and Guan \cite{Bian-Guan2009}, which are needed in the next section for the proof of the differential inequality corresponding to the minimal rank 2 cases. The proofs are omitted here and can be found in their original work.
\begin{lemma}[Lemma 2.5 in \cite{Bian-Guan2009}]\label{Lemma-v_{ijalpha}}
Assume $v\in C^{3,1}(\Omega)$ is a convex function. Then there exists a positive constant $C$ depending only on $\mathrm{dist}(\mathcal{O},\partial\Omega)$ and $\|v\|_{C^{3,1}(\Omega)}$ such that
\begin{equation}
|v_{ij\alpha}|\leq C\bigg(\sqrt{v_{ii}(x)}+\sqrt{v_{jj}(x)}\bigg)
\end{equation}
for all $x\in\mathcal{O}$ and $1\leq i,j,\alpha\leq n$.
\end{lemma}
\begin{lemma}[Lemma 3.3 in \cite{Bian-Guan2009}]\label{Lemma-v_{ijk}}
Let $M$ be a positive constant such that $0<\lambda_i\leq M$ and $\frac{1}{M}\leq\gamma_i\leq M$ for $i=l+1, l+2, \cdots, n$. Assume that $v_{ij\alpha}=v_{ji\alpha}$ for $i, j=l+1, l+2, \cdots, n$ and $\alpha=1, 2, \cdots, n$. Then there exists a positive constant $C$ depending only on $n$ and $M$, such that for each $\alpha$, for any $D>0$ and any $\delta>0$, we have
\begin{align*}
\sum_{l+1\leq i, j\leq n}|v_{ij\alpha}|
\leq&~C\left(1+\frac{2D}{\delta}+D\right)\left(\sigma_1(\lambda) +\left|\sum_{l+1\leq i\leq n}\gamma_iv_{ii\alpha}\right|\right)\\
&+\frac{\delta}{2D}\frac{1}{\sigma_1(\lambda)}\sum_{\substack{l+1\leq i, j\leq n\\ i\neq j}}|v_{ij\alpha}|^2+\frac{C}{D}\frac{1}{\sigma_1^3(\lambda)}
\sum_{l+1\leq i\leq n}V_{i\alpha}^2,
\end{align*}
where $$\sigma_1(\lambda)=\sum_{l+1\leq i\leq n}\lambda_i \quad \mathrm{and} \ V_{i\alpha}=v_{ii\alpha}\sigma_1(\lambda)-\lambda_i\bigg(\sum_{l+1\leq j\leq n}v_{jj\alpha}\bigg). $$
\end{lemma}

\section{Constant rank theorem}\label{Section 3}

This section is devoted to establishing the constant rank theorem for solutions of equation \eqref{Equation-real equation(v)}, which serves as a key ingredient in the proof of Theorem \ref{Theorem-Main}. For clarity, we adopt the convention that all summation indices $i$, $j$, $k$, $l$, $p$, $q$, $r$, $s$ and $\alpha$, $\beta$, $\gamma$, $\delta$ range from $1$ to $4$, unless stated otherwise.

\subsection{The minimal rank is $2$}\label{Section 3.1} In this subsection, we consider the case that $\nabla^2v(x)$ attains its minimal rank $2$, which is equivalent to the case where the auxiliary matrix $K$ has rank $0$. Inspired by the work of Bian and Guan \cite{Bian-Guan2009}, we employ the refined quotient type auxiliary function defined in \eqref{Auxiliary function-l=2} and derive a suitable differential inequality for $\sum_{i,j}F^{ij}\phi_{ij}$, where $F^{ij}$ denotes the linearized operator of \eqref{Equation-Dirichlet problem-real}. This estimate enables us to apply the strong maximum principle to $\phi$ and complete the proof of the constant rank theorem in this case.

\begin{theorem}\label{Theorem-minimal rank is 2}
Suppose $\Omega\subset\mathbb{R}^4$ is a domain and $v\in C^4(\Omega)$ is a convex solution of equation \eqref{Equation-real equation(v)}. Let $K(x)$ be the matrix defined by \eqref{Auxiliary matrix K}. For each $x\in\Omega$, set $$\phi(x)=\mathrm{tr}(K(x))+q(K(x)),$$ where $$q(K(x))=\begin{cases}
\frac{\det(K(x))}{\mathrm{tr}(K(x))} & \mathrm{if}\ \mathrm{tr}(K(x))>0\\
0 & \mathrm{if}\ \mathrm{tr}(K(x))=0.
\end{cases}$$ If the Hessian $\nabla^2v(x)$ attains the minimal rank $2$ at some point $x_0\in\Omega$, then there exists a neighborhood $\mathcal{O}$ of $x_0$ and positive constant $C$ independent of $\phi$, such that
\begin{equation}\label{Inequ-Theorem-minimal rank=2}
\sum_{i,j}F^{ij}\phi_{ij}\leq C(\phi+|\nabla\phi|)\quad \mathrm{in}\ \mathcal{O}.
\end{equation}
\end{theorem}

\begin{proof}
Before deriving the differential inequality for the auxiliary function, we first introduce a perturbation to remove the possible degeneracy caused by the quotient term. Following the perturbation procedure of Bian and Guan \cite{Bian-Guan2009}, for a sufficiently small $\varepsilon>0$, we define $v_\varepsilon(x)=v(x)+\frac{\varepsilon}{2}|x|^2$. Then $\nabla^2v_\varepsilon=\nabla^2v+\varepsilon I_4$, where $I_4$ denotes the $4\times 4$ identity matrix. Let $A_\varepsilon$, $B_\varepsilon$, and $K_\varepsilon$ denote the matrices associated with $v_\varepsilon$ according to the definitions of $A$, $B$, and $K$, respectively. Correspondingly, we define $$\psi_\varepsilon=\mathrm{tr}\big((\det(A_\varepsilon))^2I-B_\varepsilon\overline{A^\ast_\varepsilon}\overline{B_\varepsilon}A^\ast_\varepsilon\big),$$ and $$\phi_\varepsilon=(\det(A_\varepsilon))^{-2}\psi_\varepsilon+q_\varepsilon(\nabla^2v_\varepsilon),$$ where $$q_\varepsilon(\nabla^2v_\varepsilon)=\frac{\sigma_4(\nabla^2v_\varepsilon)}{\psi_\varepsilon}.$$ The role of this perturbation is to ensure that the denominator $\psi_\varepsilon$ is strictly positive. In fact, after choosing the coordinate system introduced below, the perturbation gives $v_{33},v_{44}\geq\varepsilon$, and hence the estimate for $\psi_\varepsilon$ (see \eqref{Estimate-psi}) yields $$\psi_\varepsilon(x)\geq c\varepsilon,\quad \mbox{for all $x\in\mathcal{O}$},$$ where $c>0$ is independent of $\varepsilon$. Since $\mathrm{tr}(K_\varepsilon)=(\det(A_\varepsilon))^{-2}\psi_\varepsilon$, it follows that
\begin{equation}\label{Estimate-varepsilon}
\mathrm{tr}(K_\varepsilon(x))\geq c\varepsilon,\quad\mbox{for all $x\in\mathcal{O}$}.
\end{equation}
Notice that $v_\varepsilon(x)$ no longer solves the original equation exactly. Instead, it satisfies
\begin{equation}\label{Equation-F-varepsilon}
F(\nabla^2v_\varepsilon(x),\nabla v_\varepsilon(x))=R_\varepsilon(x),
\end{equation}
where $R_\varepsilon(x)=F(\nabla^2v_\varepsilon,\nabla v_\varepsilon)-F(\nabla^2v,\nabla v)$ and $F(\nabla^2v,\nabla v)$ is defined as in \eqref{Equation-Dirichlet problem-real}. Since $v\in C^4(\Omega)$, for any $x\in\mathcal{O}$, we have
\begin{equation}\label{Estimate-R}
|R_\varepsilon(x)|\leq C\varepsilon,\quad |\nabla R_\varepsilon(x)|\leq C\varepsilon\ \ \mathrm{and}\ 
\ |\nabla^2R_\varepsilon(x)|\leq C\varepsilon,
\end{equation}
where $C$ is independent of $\varepsilon$.

In the following argument, we work with the perturbed equation \eqref{Equation-F-varepsilon} and the auxiliary function $\phi_\varepsilon$. For simplicity of notation, we omit the subscript $\varepsilon$ and write $v$, $\nabla^2v$, $R$, $\phi$, $A$, $q$, $\psi$ in place of $v_\varepsilon$, $\nabla^2v_\varepsilon$, $R_\varepsilon$, $\phi_\varepsilon$, $A_\varepsilon$, $q_\varepsilon$, $\psi_\varepsilon$, respectively. Under this convention, $v$ satisfies equation
\begin{equation}\label{Equation-minimal rank is 2}
F(\nabla^2v(x),\nabla v(x))=R(x).
\end{equation}
All estimates below are uniform with respect to $\varepsilon$. After deriving the required differential inequality \eqref{Inequ-Theorem-minimal rank=2} for $v_\varepsilon$, we let $\varepsilon\rightarrow 0$ to recover the corresponding result for the original function $v$.

We fix $x\in\mathcal{O}$ and choose the coordinates as in Section \ref{Section 2}. At the point $x$, the matrix $\nabla^2v$ is diagonal and satisfies $v_{11}+v_{33}=v_{22}+v_{44}$, with $v_{11}(x)\geq v_{22}(x)\geq v_{33}(x)\geq v_{44}(x)\geq 0$. It follows from \eqref{Estimate-varepsilon} that $\varepsilon\leq C\phi(x)$. Hence, using the estimates in \eqref{Estimate-R}, we obtain
\begin{equation}\label{Estimate-R1}
R(x)=O(\phi(x)),\quad \nabla R(x)=O(\phi(x))\ \ \mathrm{and}\ \ \nabla^2R(x)=O(\phi(x)).
\end{equation}
Throughout the proof, we perform all calculations at the fixed point $x$. Since we are considering the case where $\nabla^2v$ has minimal rank $2$, we use the notation $G=\{1,2\}$, $B=\{3,4\}$ for the “good" and “bad" indices, respectively. The good eigenvalues are uniformly positive; namely, $v_{ii}\geq c$ for $i\in G$. For convenience, we also restate the previously established estimates \eqref{Estimate-v_{33},v_{44}} and \eqref{Estimate-v_{22}} as
\begin{equation}\label{Estimate-v_{33},v_{44}'}
v_{33}=O(\phi),\quad v_{44}=O(\phi)
\end{equation}
and
\begin{equation}\label{Estimate-v_{22}'}
v_{22}=v_{11}+O(\phi).
\end{equation}

Taking the first derivatives of $\phi$ at the fixed point $x$, and applying Lemma \ref{Lemma-first derivative of q}, we obtain
\begin{align*}
\phi_i
=&-2\cdot16^2(v_{11}+v_{33})^{-3}(v_{22}+v_{44})^{-2}(v_{11i}+v_{33i})\psi\\
&-2\cdot16^2(v_{11}+v_{33})^{-2}(v_{22}+v_{44})^{-3}(v_{22i}+v_{44i})\psi\\
&+16^2(v_{11}+v_{33})^{-2}(v_{22}+v_{44})^{-2}\psi_i+O\bigg(\phi+\sum_{\alpha,\beta\in B}|\nabla v_{\alpha\beta}|\bigg).
\end{align*}
For $\psi_i$, recalling the expression of $\psi$ in \eqref{Equality-psi}, a direct computation gives
\begin{align*}
\psi_i
=&~\frac{4}{16^2}(v_{11}+v_{33})^2(v_{22}+v_{44})(v_{22i}+v_{44i})+\frac{4}{16^2}(v_{11}+v_{33})(v_{22}+v_{44})^2(v_{11i}+v_{33i})\\
&-\frac{2}{16^2}(v_{11}-v_{33})(v_{22}+v_{44})^2(v_{11i}-v_{33i})-\frac{2}{16^2}(v_{11}-v_{33})^2(v_{22}+v_{44})(v_{22i}+v_{44i})\\
&-\frac{2}{16^2}(v_{22}-v_{44})(v_{11}+v_{33})^2(v_{22i}-v_{44i})-\frac{2}{16^2}(v_{22}-v_{44})^2(v_{11}+v_{33})(v_{11i}+v_{33i}).
\end{align*}
By means of the estimates \eqref{Estimate-v_{33},v_{44}'} and \eqref{Estimate-v_{22}'}, this implies
\begin{equation}\label{Equality-psi_i}
\begin{aligned}
\psi_i
=&~\frac{4}{16^2}v_{11}^3(v_{22i}+v_{44i})+\frac{4}{16^2}v_{11}^3(v_{11i}+v_{33i})-\frac{2}{16^2}v_{11}^3(v_{11i}-v_{33i})\\
&-\frac{2}{16^2}v_{11}^3(v_{22i}+v_{44i})-\frac{2}{16^2}v_{11}^3(v_{22i}-v_{44i})-\frac{2}{16^2}v_{11}^3(v_{11i}+v_{33i})+O(\phi)\\
=&~\frac{4}{16^2}v_{11}^3(v_{33i}+v_{44i})+O(\phi)=O\bigg(\phi+\sum_{p,q\in B}|\nabla v_{pq}|\bigg).
\end{aligned}
\end{equation}
Moreover, since $v_{11}+v_{33}$ and $v_{22}+v_{44}$ are uniformly bounded and $\psi=O(\phi)$, we conclude that
\begin{equation}\label{Equality-phi_i-2}
\phi_i=O\bigg(\phi+\sum_{\alpha,\beta\in B}|\nabla v_{\alpha\beta}|\bigg).
\end{equation}

Next, differentiating $\phi$ twice and evaluating at $x$, we obtain
\begin{align*}
\phi_{ij}
=&~6\cdot16^2(v_{11}+v_{33})^{-4}(v_{22}+v_{44})^{-4}\\
&\cdot\big[(v_{11i}+v_{33i})(v_{22}+v_{44})+(v_{11}+v_{33})(v_{22i}+v_{44i})\big]\\
&\cdot\big[(v_{11j}+v_{33j})(v_{22}+v_{44})+(v_{11}+v_{33})(v_{22j}+v_{33j})\big]\psi\\
&-2\cdot16^2(v_{11}+v_{33})^{-3}(v_{22}+v_{44})^{-3}\\
&\cdot\big[(v_{11ij}+v_{33ij})(v_{22}+v_{44})+(v_{11}+v_{33})(v_{22ij}+v_{44ij})\\
&\hspace{0.4cm}+(v_{11i}v_{33i})(v_{22j}+v_{44j})+(v_{11j}+v_{33j})(v_{22i}+v_{44i})\\
&\hspace{0.4cm}-2(v_{12i}+v_{34i})(v_{12j}+v_{34j})-2(v_{14i}-v_{32i})(v_{14j}-v_{32j})\big]\psi\\
&-2\cdot16^2(v_{11}+v_{33})^{-3}(v_{22}+v_{44})^{-3}\\
&\cdot\big[(v_{11i}+v_{33i})(v_{22}+v_{44})+(v_{11}+v_{33})(v_{22i}+v_{44i})\big]\psi_j\\
&-2\cdot16^2(v_{11}+v_{33})^{-3}(v_{22}+v_{44})^{-3}\\
&\cdot\big[(v_{11j}+v_{33j})(v_{22}+v_{44})+(v_{11}+v_{33})(v_{22j}+v_{44j})\big]\psi_i\\
&+16^2(v_{11}+v_{33})^{-2}(v_{22}+v_{44})^{-2}\psi_{ij}+q_{ij}.
\end{align*}
The first two terms on the right-hand side are bounded by $O(\phi)$, due to the uniform bounds of $v_{11}+v_{33}$ and $v_{22}+v_{44}$, together with $\psi=O(\phi)$. The third and fourth terms are estimated by $O(\phi+\sum_{\alpha,\beta\in B}|\nabla v_{\alpha\beta}|)$ in view of \eqref{Equality-psi_i}. For the fifth terms, using \eqref{Estimate-v_{33},v_{44}'} and \eqref{Estimate-v_{22}'}, we have $$16^2(v_{11}+v_{33})^{-2}(v_{22}+v_{44})^{-2}=16^2v_{11}^{-4}+O(\phi).$$ It remains to calculate $\psi_{ij}$. Since $\nabla^2v$ is diagonal at the point under consideration, direct differentiation of the expression of $\psi$ gives
\begin{align*}
\psi_{ij}
=&~\frac{4}{16^2}\big[(v_{11i}+v_{33i})(v_{22}+v_{44})+(v_{11}+v_{33})(v_{22i}+v_{44i})\big]\\
&\hspace{0.8cm}\cdot\big[(v_{11j}+v_{33j})(v_{22}+v_{44})+(v_{11}+v_{33})(v_{22j}+v_{44j})\big]\\
&+\frac{4}{16^2}(v_{11}+v_{33})(v_{22}+v_{44})\\
&\hspace{0.8cm}\cdot\big[(v_{11ij}+v_{33ij})(v_{22}+v_{44})+(v_{11}+v_{33})(v_{22ij}+v_{44ij})\\
&\hspace{1.2cm}+(v_{11i}+v_{33i})(v_{22j}+v_{44j})+(v_{11j}+v_{33j})(v_{22i}+v_{44\i})\\
&\hspace{1.2cm}-2(v_{12i}+v_{34i})(v_{12j}+v_{34j})-2(v_{14i}-v_{32i})(v_{14j}-v_{32j})\big]\\
&-\frac{1}{16^2}\big\{2\big[(v_{11i}-v_{33i})(v_{22}+v_{44})+(v_{11}-v_{33})(v_{22i}+v_{44i})\big]\\
&\hspace{1.3cm}\cdot\big[(v_{11j}-v_{33j})(v_{22}+v_{44})+(v_{11}-v_{33})(v_{22j}+v_{44j})\big]\\
&\hspace{1.2cm}+2(v_{11}-v_{33})(v_{22}+v_{44})\\
&\hspace{1.4cm}\cdot\big[(v_{11ij}-v_{33ij})(v_{22}+v_{44})+(v_{11}-v_{33})(v_{22ij}+v_{44ij})\\
&\hspace{1.7cm}+(v_{11i}-v_{33i})(v_{22j}+v_{44j})+(v_{11j}-v_{33j})(v_{22i}+v_{44i})\\
&\hspace{1.7cm}+2v_{23i}v_{23j}+2v_{34i}v_{34j}-2v_{12i}v_{12j}-2v_{14i}v_{14j}\big]\\
&\hspace{1.2cm}+2\big[(v_{22i}-v_{44i})(v_{11}+v_{33})+(v_{22}-v_{44})(v_{11i}+v_{33i})\big]\\
&\hspace{1.6cm}\cdot\big[(v_{22j}-v_{44j})(v_{11}+v_{33})+(v_{22}-v_{44})(v_{11j}+v_{33j})\big]\\
&\hspace{1.2cm}+2(v_{22}-v_{44})(v_{11}+v_{33})\\
&\hspace{1.6cm}\cdot\big[(v_{22ij}-v_{44ij})(v_{11}+v_{33})+(v_{22}-v_{44})(v_{11ij}+v_{33ij})\\
&\hspace{1.8cm}+(v_{22i}-v_{44i})(v_{11j}+v_{33j})+(v_{22j}-v_{44j})(v_{11i}+v_{33i})\\
&\hspace{1.8cm}+2v_{14i}v_{14j}+2v_{34i}v_{34j}-2v_{12i}v_{12j}-2v_{23i}v_{23j}\big]\\
&\hspace{1.2cm}+8(v_{22}+v_{44})^2v_{12i}v_{13j}+8(v_{11}+v_{33})^2v_{24i}v_{24j}\\
&\hspace{1.2cm}+8(v_{12i}v_{33}-v_{11}v_{34i})(v_{21j}v_{44}-v_{22}v_{43j})\\
&\hspace{1.2cm}+8(v_{12j}v_{33}-v_{11}v_{34j})(v_{21i}v_{44}-v_{22}v_{43i})\\
&\hspace{1.2cm}+8(v_{11}v_{32i}+v_{33}v_{14i})(v_{44}v_{23j}+v_{22}v_{14j})\\
&\hspace{1.2cm}+8(v_{11}v_{32j}+v_{33}v_{14j})(v_{44}v_{23i}+v_{22}v_{14i})\big\}.
\end{align*}
Applying Lemma \ref{Lemma-v_{ijalpha}} and the Cauchy-Schwarz inequality, we derive the estimate $$v_{pq\alpha}^2=O(\phi),\quad p,q\in B,\ \alpha=\{1,2,3,4\}.$$ Combining this estimate with \eqref{Estimate-v_{33},v_{44}'} and \eqref{Estimate-v_{22}'}, we further obtain the following estimate for $\psi_{ij}$
\begin{align*}
\psi_{ij}
=&~\frac{4}{16^2}v_{11}^3(v_{33ij}+v_{44ij})-\frac{8}{16^2}v_{11}^2v_{13i}v_{13j}-\frac{8}{16^2}v_{11}^2v_{14i}v_{14j}\\
&-\frac{8}{16^2}v_{11}^2v_{32i}v_{32j}-\frac{8}{16^2}v_{11}^2v_{24i}v_{24j}+O\bigg(\phi+\sum_{p,q\in B}|\nabla v_{pq}|\bigg).
\end{align*}
Therefore, we conclude that
\begin{align*}
16^2(v_{11}+v_{33})^{-2}(v_{22}+v_{44})^{-2}\psi_{ij}=&~4v_{11}^{-1}(v_{33ij}+v_{44ij})-8v_{11}^{-2}v_{13i}v_{13j}-8v_{11}^{-2}v_{14i}v_{14j}\\
&-8v_{11}^{-2}v_{32i}v_{32j}-8v_{11}^{-2}v_{24i}v_{24j}+O\bigg(\phi+\sum_{\alpha,\beta\in B}|\nabla v_{\alpha\beta}|\bigg).
\end{align*}
Collecting the above estimates and applying the formula for $q_{ij}$ established in Lemma \ref{Lemma-q_{ij}}, we finally obtain
\begin{align*}
\phi_{ij}
=&~4v_{11}^{-1}\left(1+16\frac{v_{44}^2}{\sigma_1^2(B)}\right)v_{33ij}+4v_{11}^{-1}\left(1+16\frac{v_{33}^2}{\sigma_1^2(B)}\right)v_{44ij}\\
&-8v_{11}^{-2}\left(1+16\frac{v_{44}^2}{\sigma_1^2(B)}\right)v_{13i}v_{13j}-8v_{11}^{-2}\left(1+16\frac{v_{33}^2}{\sigma_1^2(B)}\right)v_{14i}v_{14j}\\
&-8v_{11}^{-2}\left(1+16\frac{v_{44}^2}{\sigma_1^2(B)}\right)v_{23i}v_{23j}-8v_{11}^{-2}\left(1+16\frac{v_{33}^2}{\sigma_1^2(B)}\right)v_{24i}v_{24j}\\
&-4\cdot16v_{11}^{-1}\frac{1}{\sigma_1^3(B)}(V_{3i}V_{3j}+V_{4i}V_{4j})-8\cdot16v_{11}^{-1}\frac{1}{\sigma_1(B)}v_{34i}v_{34j}\\
&+O\bigg(\phi+\sum_{\alpha,\beta\in B}|\nabla v_{\alpha\beta}|\bigg),
\end{align*}
where $$V_{3i}=v_{33i}\sigma_1(B)-v_{33}\bigg(\sum_{\gamma\in B}v_{\gamma\gamma i}\bigg),\quad\mathrm{and}\ V_{4i}=v_{44i}\sigma_1(B)-v_{44}\bigg(\sum_{\gamma\in B}v_{\gamma\gamma i}\bigg).$$ Contracting this identity with $F^{ij}$, we arrive at
\begin{equation}\label{Equality-F1-2}
\begin{aligned}
\sum_{i,j}F^{ij}\phi_{ij}
=&~4v_{11}^{-1}\sum_{i,j}\sum_{\alpha\in B}\left(1+16\frac{\sigma_1^2(B|\alpha)}{\sigma_1^2(B)}\right)F^{ij}v_{\alpha\alpha ij}\\
&-8v_{11}^{-2}\sum_{i,j}\sum_{\alpha\in B,\beta\in G}\left(1+16\frac{\sigma_1^2(B|\alpha)}{\sigma_1^2(B)}\right)F^{ij}v_{\alpha\beta i}v_{\alpha\beta j}\\
&-4\cdot 16v_{11}^{-1}\frac{1}{\sigma_1^3(B)}\sum_{i,j}\sum_{\alpha\in B}F^{ij}V_{\alpha i}V_{\alpha j}\\
&-4\cdot16v_{11}^{-1}\frac{1}{\sigma_1(B)}\sum_{i,j}\sum_{\substack{\alpha,\beta\in B\\ \alpha\neq\beta}}F^{ij}v_{\alpha\beta i}v_{\alpha\beta j}\\
&+O\bigg(\phi+\sum_{\alpha,\beta\in B}|\nabla v_{\alpha\beta}|\bigg),
\end{aligned}
\end{equation}
where $$V_{\alpha i}=v_{\alpha\alpha i}\sigma_1(B) -v_{\alpha\alpha}\bigg( \sum_{\gamma\in B}v_{\gamma\gamma i}\bigg).$$
 
To manage the fourth-order derivative terms in \eqref{Equality-F1-2}, we differentiate equation \eqref{Equation-minimal rank is 2} twice to obtain
\begin{equation}\label{Relation-1st condition-2}
\sum_{i, j}F^{ij}v_{ij\alpha}+\sum_kF^{v_k}v_{k\alpha}=O(\phi),
\end{equation}
and
\begin{align*}
O(\phi)=&\sum_{i,j,k,l}F^{ij,kl}v_{ij\alpha}v_{kl\beta}+\sum_{i,j,k}F^{ij,v_k}v_{ij\alpha}v_{k\beta}+\sum_{i,j}F^{ij}v_{ij\alpha\beta}\\
&+\sum_{k,i,j}F^{v_k,ij}v_{k\alpha}v_{ij\beta}+\sum_{k,l}F^{v_k,v_l}v_{k\alpha}v_{l\beta}+\sum_kF^{v_k}v_{k\alpha\beta},
\end{align*}
where we used the estimates in \eqref{Estimate-R1}. Therefore, it follows from the estimates \eqref{Estimate-v_{33},v_{44}'}, \eqref{Estimate-v_{22}'}, and Lemma \ref{Lemma-v_{ijalpha}} that
\begin{equation}\label{Equality-4th order derivative}
\begin{aligned}
&~4v_{11}^{-1}\sum_{i,j}\sum_{\alpha\in B}\left(1+16\frac{\sigma_1^2(B|\alpha)}{\sigma_1^2(B)}\right)F^{ij}v_{\alpha\alpha ij}\\
=&-4v_{11}^{-1}\sum_{\alpha\in B}\left(1+16\frac{\sigma_1^2(B|\alpha)}{\sigma_1^2(B)}\right)\\
&\cdot\bigg(\sum_{i,j,k,l}F^{ij,kl}v_{ij\alpha}v_{kl\alpha}+2\sum_{i,j,k}F^{ij,v_k}v_{ij\alpha}v_{k\alpha}+\sum_{k,l}F^{v_k,v_l}v_{k\alpha}v_{l\alpha}+\sum_kF^{v_k}v_{k\alpha\alpha}+O(\phi)\bigg)\\
=&-4v_{11}^{-1}\sum_{\alpha\in B}\left(1+16\frac{\sigma_1^2(B|\alpha)}{\sigma_1^2(B)}\right)\bigg(\sum_{i,j,k,l\in G}F^{ij,kl}v_{ij\alpha}v_{kl\alpha}\bigg)+O\bigg(\phi+\sum_{\alpha,\beta\in B}|\nabla v_{\alpha\beta}|\bigg).
\end{aligned}
\end{equation}
Substituting \eqref{Equality-4th order derivative} into \eqref{Equality-F1-2} yields
\begin{equation}\label{Equality-F2-2}
\begin{aligned}
\sum_{i,j}F^{ij}\phi_{ij}
=&~-4v_{11}^{-1}\sum_{\alpha\in B}\left(1+16\frac{\sigma_1^2(B|\alpha)}{\sigma_1^2(B)}\right)\\
&\cdot\bigg(\sum_{i,j,k,l\in G}F^{ij,kl}v_{ij\alpha}v_{kl\alpha}+2v_{11}^{-1}\sum_{i,j,\beta\in G}F^{ij}v_{\alpha\beta i}v_{\alpha\beta j}\bigg)\\
&-4\cdot16v_{11}^{-1}\bigg(\frac{1}{\sigma_1^3(B)}\sum_{i,j}\sum_{\alpha\in B}F^{ij}V_{\alpha i}V_{\alpha j}+\frac{1}{\sigma_1(B)}\sum_{i,j}\sum_{\substack{\alpha,\beta\in B\\ \alpha\neq\beta}}F^{ij}v_{\alpha\beta i}v_{\alpha\beta j}\bigg)\\
&+O\bigg(\phi+\sum_{\alpha,\beta\in B}|\nabla v_{\alpha\beta}|\bigg).
\end{aligned}
\end{equation}
 
In order to obtain the desired differential inequality \eqref{Inequ-Theorem-minimal rank=2}, we shall prove the following claim.
 
\textbf{Claim 1.} For each $\alpha\in B$, we have
\begin{equation}\label{Claim}
\sum_{i,j,k,l\in G}F^{ij,kl}v_{ij\alpha}v_{kl\alpha}+2\sum_{i,j,\beta\in G}\frac{1}{v_{\beta\beta}}F^{ij}v_{\alpha\beta i}v_{\beta\alpha j}\geq -C\bigg(\phi+\sum_{\alpha,\beta\in B}|\nabla v_{\alpha\beta}|\bigg).
\end{equation}
 
Once \textbf{Claim 1} is proved, since $v_{22}=v_{11}+O(\phi)$, \eqref{Equality-F2-2} implies that
\begin{align*}
\sum_{i, j}F^{ij}\phi_{ij} \leq&-64v_{11}^{-1}\bigg(\frac{1}{\sigma_1^3(B)}\sum_{i, j} \sum_{\alpha\in B}F^{ij}V_{\alpha i}V_{\alpha j} +\frac{1}{\sigma_1(B)} \sum_{i,j} \sum_{\substack{\alpha, \beta\in B\\ \alpha\neq\beta}} F^{ij}v_{\alpha\beta i}v_{\alpha\beta j}\bigg)\\
&+C\bigg(\phi+\sum_{\alpha, \beta\in B}|\nabla v_{\alpha\beta}|\bigg),
\end{align*}
where $V_{\alpha i}=v_{\alpha\alpha i}\sigma_1(B)-v_{\alpha\alpha} \left(\sum_{\gamma\in B}v_{\gamma\gamma i}\right)$. Since $v\in C^4(\Omega)$ and $\bar{\mathcal{O}}\subset\Omega$, there exists a constant $\delta_0>0$ such that $$(F^{ij})\geq\delta_0I_4\quad \mbox{in}\ \mathcal{O}. $$ Consequently, we can deduce $$\sum_{i,j}F^{ij}V_{\alpha i}V_{\alpha j}\geq\delta_0\sum_i V_{\alpha i}^2 \quad \mathrm{and} \quad \sum_{i,j}F^{ij} v_{\alpha\beta i}v_{\beta\alpha j}\geq\delta_0\sum_i v_{\alpha\beta i}^2. $$ It follows that $$\sum_{i,j}F^{ij}\phi_{ij}\leq -\frac{\delta_0'}{\sigma_1^3(B)}\sum_i\sum_{\alpha\in B}V_{\alpha i}^2-\frac{\delta_0'}{\sigma_1(B)} \sum_i\sum_{\substack{\alpha, \beta\in B\\ \alpha\neq\beta}}v_{\alpha\beta i}^2+C\bigg(\phi+\sum_{\alpha,\beta\in B}|\nabla v_{\alpha\beta}|\bigg). $$ By Lemma \ref{Lemma-v_{ijk}} and equation \eqref{Equality-phi_i-2}, we finally obtain $$\sum_{i, j}F^{ij}\phi_{ij}\leq C(\phi+|\nabla\phi|). $$ Since the constant $C$ is independent of $\varepsilon$, letting $\varepsilon\rightarrow 0$ yields \eqref{Inequ-Theorem-minimal rank=2} for $v$. 
\end{proof}

\textbf{Proof of Claim.} Recall the definition of $F$, namely,
\begin{align*}
F(\nabla^2v, \nabla v)=&~(v_{11}+v_{33})(v_{22}+v_{44}) -(v_{12}+v_{34})(v_{21}+v_{43}) -(v_{14}-v_{32})(v_{41}-v_{23})\\
&-(v_2^2+v_4^2)(v_{11}+v_{33})-(v_1^2+v_3^2)(v_{22}+v_{44})\\
&+(v_1v_2+v_3v_4)(v_{12}+v_{34}+v_{21}+v_{43})\\
&+(v_1v_4-v_2v_3)(v_{14}-v_{32}+v_{41}-v_{23})-\lambda.
\end{align*}
It is clear that $$F^{11,22}=F^{22,11}=1, \quad F^{12,21}=F^{21,12}=-1. $$ Note that $G=\{1, 2\}$ and $B=\{3, 4\}$. Expanding and collecting the terms on the left-hand side of \eqref{Claim}, we obtain
\begin{equation}\label{Equality-F3-2}
\begin{aligned}
&\sum_{i,j,k,l\in G}F^{ij,kl}v_{ij\alpha}v_{kl\alpha}+2\sum_{i,j,\beta\in G}\frac{1}{v_{\beta\beta}}F^{ij}v_{\alpha\beta i}v_{\beta\alpha j}\\
=&~\frac{2F^{11}}{v_{11}}v_{11\alpha}^2+\frac{2}{v_{11}v_{22}}\left(F^{11}v_{11}+F^{22}v_{22}-v_{11}v_{22}\right)v_{12\alpha}^2+\frac{2F^{22}}{v_{22}}v_{22\alpha}^2\\
&+\frac{4F^{12}}{v_{11}}v_{11\alpha}v_{12\alpha}+2v_{11\alpha}v_{22\alpha}+\frac{4F^{12}}{v_{22}}v_{12\alpha}v_{22\alpha}.
\end{aligned}
\end{equation}
By \eqref{Estimate-v_{33},v_{44}'}, equation \eqref{Relation-1st condition-2} is equivalent to $$\sum_{i,j\in G}F^{ij}v_{ij\alpha}=O\bigg(\phi +\sum_{\alpha,\beta\in B}|\nabla v_{\alpha\beta}|\bigg), $$ namely, 
\begin{equation}\label{Relation-v_22alpha-2}
v_{22\alpha}=-\frac{F^{11}}{F^{22}}v_{11\alpha}-2\frac{F^{12}}{F^{22}}v_{12\alpha}+O\bigg(\phi+\sum_{\alpha,\beta\in B}|\nabla v_{\alpha\beta}|\bigg).
\end{equation}
In equation \eqref{Relation-v_22alpha-2}, we used the fact that $F^{22}\neq0$. This is ensured by the ellipticity of equation \eqref{Equation-Dirichlet problem-real}. Putting \eqref{Relation-v_22alpha-2} into \eqref{Equality-F3-2}, we have
\begin{equation}\label{Equality-claim-2}
\begin{aligned}
&\sum_{i,j,k,l\in G}F^{ij,kl}v_{ij\alpha}v_{kl\alpha}+2\sum_{i,j,\beta\in G}\frac{1}{v_{\beta\beta}}F^{ij}v_{\alpha\beta i}v_{\beta\alpha j}\\
=&~\frac{2F^{11}}{v_{11}}\frac{1}{F^{22}v_{22}}\left(F^{11}v_{11}+F^{22}v_{22}-v_{11}v_{22}\right)v_{11\alpha}^2\\
&+\frac{4F^{12}}{v_{11}}\frac{1}{F^{22}v_{22}}\left(F^{11}v_{11}+F^{22}v_{22}-v_{11}v_{22}\right)v_{11\alpha}v_{12\alpha}\\
&+\frac{2}{v_{11}v_{22}}\left(F^{11}v_{11}+F^{22}v_{22}-v_{11}v_{22}\right)v_{12\alpha}^2+O\bigg(\phi+\sum_{\alpha,\beta\in B}|\nabla v_{\alpha\beta}|\bigg).
\end{aligned}
\end{equation}
 
Using equation \eqref{Equation-minimal rank is 2} along with the estimates $v_{33}=O(\phi)$ and $v_{44}=O(\phi)$, we readily obtain $$v_{11}v_{22}-(v_1^2+v_3^2)v_{22}-(v_2^2+v_4^2)v_{11}=\lambda+O(\phi), $$ and $$F^{11}=v_{22}-(v_2^2+v_4^2)+O(\phi), \quad F^{22}=v_{11}-(v_1^2+v_3^2)+O(\phi). $$ Consequently, $$F^{11}v_{11}+F^{22}v_{22}-v_{11}v_{22}=\lambda+O(\phi). $$ Substituting this into \eqref{Equality-claim-2}, we obtain
\begin{align}\label{Equality-claim-3}
\begin{split}
&\sum_{i,j,k,l\in G}F^{ij,kl}v_{ij\alpha}v_{kl\alpha}+2\sum_{i,j,\beta\in G}\frac{1}{v_{\beta\beta}}F^{ij}v_{\alpha\beta i}v_{\beta\alpha j}\\
=&~\frac{2\lambda}{v_{11}\cdot F^{22}v_{22}}\left(F^{11}v_{11\alpha}^2+2F^{12}v_{11\alpha}v_{12\alpha}+F^{22}v_{12\alpha}^2\right)+O\bigg(\phi+\sum_{\alpha,\beta\in B}|\nabla v_{\alpha\beta}|\bigg)\\
=&~\frac{2\lambda}{v_{11}\cdot F^{22}v_{22}}\left(\sqrt{F^{11}}v_{11\alpha}+\frac{F^{12}}{\sqrt{F^{11}}}v_{12\alpha}\right)^2\\
&+\frac{2\lambda}{F^{11}v_{11}\cdot F^{22}v_{22}}\left(F^{11}F^{22}-(F^{12})^2\right)v_{12\alpha}^2+O\bigg(\phi+\sum_{\alpha,\beta\in B}|\nabla v_{\alpha\beta}|\bigg).
\end{split}
\end{align}
Since $F^{12}=v_1v_2+v_3v_4$ and $F^{23}=-(v_1v_4-v_2v_3)$, it follows that $$(F^{12})^2+(F^{23})^2=(v_1^2+v_3^2)(v_2^2+v_4^2). $$ A straightforward computation yields
\begin{align}\label{Equ:F11F22}
\begin{split}
F^{11}F^{22}=&~v_{11}v_{22}-(v_1^2+v_3^2)v_{22}-(v_2^2+v_4^2)v_{11}+(v_1^2+v_3^2)(v_2^2+v_4^2)+O(\phi)\\
=&~\lambda+(F^{12})^2+(F^{23})^2+O(\phi).
\end{split}
\end{align}
By \eqref{Equality-claim-3} and \eqref{Equ:F11F22}, we conclude that
\begin{align*}
&\sum_{i,j,k,l\in G}F^{ij,kl}v_{ij\alpha}v_{kl\alpha}+2\sum_{i,j,\beta\in G}\frac{1}{v_{\beta\beta}}F^{ij}v_{\alpha\beta i}v_{\beta\alpha j}\\
=&~\frac{2\lambda}{v_{11}\cdot F^{22}v_{22}}\left(F^{11}v_{11\alpha}^2+2F^{12}v_{11\alpha}v_{12\alpha}+F^{22}v_{12\alpha}^2\right)+O\bigg(\phi+\sum_{\alpha,\beta\in B}|\nabla v_{\alpha\beta}|\bigg)\\
=&~\frac{2\lambda}{v_{11}\cdot F^{22}v_{22}}\left(\sqrt{F^{11}}v_{11\alpha}+\frac{F^{12}}{\sqrt{F^{11}}}v_{12\alpha}\right)^2\\
&+\frac{2\lambda}{F^{11}v_{11}\cdot F^{22}v_{22}}\left(\lambda+(F^{23})^2\right)v_{12\alpha}^2+O\bigg(\phi+\sum_{\alpha,\beta\in B}|\nabla v_{\alpha\beta}|\bigg).
\end{align*}
This completes the proof of \textbf{Claim 1}. \qed

\subsection{The minimal rank is $3$}\label{Section 3.2}
In this subsection, we focus on the case that $\nabla^2v(x)$ attains its minimal rank $3$ in $\Omega$. Equivalently, the auxiliary matrix $K$ has minimal rank $1$. For this case, we consider the auxiliary function $$\phi(x)=\det(K(x)). $$ The main purpose is to derive the differential inequality \eqref{Inequ-Theorem-minimal rank=3} for $\phi(x)$ and then apply it to prove the corresponding constant rank theorem.
\begin{theorem}\label{Theorem-minimal rank is 3}
Let $\Omega\subset\mathbb{R}^4$ be a domain and let $v$ be a convex solution to equation \eqref{Equation-real equation(v)}. If the Hessian $\nabla^2v(x)$ attains its minimal rank $3$ at some point $x_0\in\Omega$, then there exist a neighborhood $\mathcal{O}$ of $x_0$ and a positive constant $C$ (independent of $\phi$) such that
\begin{equation}\label{Inequ-Theorem-minimal rank=3}
\sum_{i, j}F^{ij}\phi_{ij}\leq C(\phi+|\nabla\phi|)\quad \mathrm{in} \ \mathcal{O}.
\end{equation}
\end{theorem}
\begin{proof}
Following the notations in \cite{Caffarelli-Friedman1985, Korevaar-Lewis1987}, for two functions $h(x)$ and $k(x)$ defined in $\mathcal{O}$, we say that $h(x)\lesssim k(x)$ if there exists a positive constant $C_1$ such that $$h(x)-k(x)\leq C_1(\phi(x)+|\nabla\phi|(x)). $$ We also write $h(x)\sim k(x)$ if $h(x)\lesssim k(x)$ and $k(x)\lesssim h(x)$. Next, we write $h\lesssim k$ if the above inequality holds for any $x\in\mathcal{O}$, with constant $C_1$ independent of $x$. Finally, $h\sim k$ means $h\lesssim k$ and $k\lesssim h$. In a similar way, we will also use the notation $h\gtrsim k$ which means $-h\lesssim -k$, namely, $$-h+k\leq C_2(\phi+|\nabla\phi|)\quad \mathrm{in}\ \mathcal{O}, $$ or equivalently, $$h-k\geq-C_2(\phi+|\nabla\phi|)\quad \mathrm{in}\ \mathcal{O}. $$ Thus, $h\sim k$ can also be interpreted as $h\lesssim k$ and $h\gtrsim k$, meaning that there exist positive constants $C_3$ and $C_4$ such that $$-C_3(\phi+|\nabla\phi|) \leq h-k \leq C_4(\phi+|\nabla\phi|) \quad \mathrm{in} \ \mathcal{O}. $$ We can also write this as $$h=k+O(\phi+|\nabla\phi|). $$
 
For each $x\in\mathcal{O}$ fixed, we choose orthonormal coordinates such that $\nabla^2v(x)$ is diagonal. In the following, all calculations will be carried out at the fixed point $x$. When we use the relation $\lesssim$ or $\gtrsim$, all constants involved are under controlled. In this subsection, we set $G=\{1, 2, 3\}$ and $B=\{4\}$. According to \eqref{Equality-det K and det D^2v}, we have $$\phi=(\det(A))^{-2}\sigma_4(\nabla^2v).$$ At the fixed point, since $\nabla^2v$ has minimal rank $3$, we have $$0\sim\phi=(\det(A))^{-2}\sigma_3(G)v_{44},$$ and hence
\begin{equation}\label{Relation-v_44-3}
v_{44}\sim 0.
\end{equation}
Taking the first derivatives of $\phi$, and applying Lemma \ref{Lemma-sigma_k}, we obtain \begin{align*}
0\sim\phi_i
=&~\big((\det(A))^{-2}\big)_i\sigma_4+(\det(A))^{-2}\sum_{p,q}\frac{\partial\sigma_4}{\partial v_{pq}}v_{pqi}\\
\sim&~(\det(A))^{-2}\sum_p\sigma_3(\lambda|p)v_{ppi}\sim(\det(A))^{-2}\sigma_3(G)v_{44i},\quad i=1,2,3,4.
\end{align*}
Namely,
\begin{equation}\label{Relation-v_44i-3}
v_{44i}\sim0, \quad i=1, 2, 3, 4.
\end{equation}
Using relations \eqref{Relation-v_44-3} and \eqref{Relation-v_44i-3}, by Lemma \ref{Lemma-sigma_k} again, we calculate the second derivatives of $\phi$ as
\begin{align*}
\phi_{ij}
=&~\big((\det(A))^{-2}\big)_{ij}\sigma_4+\big((\det(A))^{-2}\big)_i\sum_{p,q}\frac{\partial\sigma_4}{\partial v_{pq}}v_{pqj}+\big((\det(A))^{-2}\big)_j\sum_{p,q}\frac{\partial\sigma_4}{\partial v_{pq}}v_{pqi}\\
&+(\det(A))^{-2}\sum_{p,q}\frac{\partial\sigma_4}{\partial v_{pq}}v_{pqij}+(\det(A))^{-2}\sum_{p,q,r,s}\frac{\partial^2\sigma_4}{\partial v_{pq}\partial v_{rs}}v_{pqi}v_{rsj}\\
=&~\big((\det(A))^{-2}\big)_{ij}\sigma_3(G)v_{44}+\big((\det(A))^{-2}\big)_i\sum_p\sigma_3(\lambda|p)v_{ppj}+\big((\det(A))^{-2}\big)_j\sum_p\sigma_3(\lambda|p)v_{ppi}\\
&+(\det(A))^{-2}\sum_p\sigma_3(\lambda|p)v_{ppij}\\
&+(\det(A))^{-2}\sum_{p\neq r}\sigma_2(\lambda|pr) v_{ppi}v_{rrj} -(\det(A))^{-2}\sum_{p\neq r}\sigma_2(\lambda|pr)v_{pri}v_{rpj}\\
\sim&~(\det(A))^{-2}\sigma_3(G)v_{44ij}-2(\det(A))^{-2}\sigma_3(G)\sum_{p\in G}\frac{1}{v_{pp}}v_{4pi}v_{4pj}
\end{align*}
It follows that
\begin{equation}\label{Equality-F1-3}
-\frac{(\det(A))^2}{\sigma_3(G)}\sum_{i,j}F^{ij}\phi_{ij}\sim-\sum_{i,j}F^{ij}v_{44ij}+2\sum_{p\in G}\sum_{i,j}\frac{1}{v_{pp}}F^{ij}v_{4pi}v_{4pj}.
\end{equation}
 
To handle the fourth-order derivative terms in \eqref{Equality-F1-3}, we first take the derivative of the equation $F(\nabla^2v, \nabla v)=0$ with respect to the variable $x_4$, resulting in
\begin{equation}\label{Equality-1st order derivative of F-3}
0=\sum_{i, j}F^{ij}v_{ij4}+\sum_kF^{v_k}v_{k4}.
\end{equation}
Then we take the second derivative, which gives
\begin{equation}\label{Equality-2nd order derivative of F-3}
\begin{aligned}
0=&\sum_{i,j,k,l}F^{ij,kl}v_{ij4}v_{kl4}+\sum_{i,j,k}F^{ij,v_k}v_{ij4}v_{k4}+\sum_{i,j}F^{ij}v_{ij44}\\
&+\sum_{k,i,j}F^{v_k,ij}v_{k4}v_{ij4}+\sum_{k,l}F^{v_k,v_l}v_{k4}v_{l4}+\sum_kF^{v_k}v_{k44}.
\end{aligned}
\end{equation}
By \eqref{Relation-v_44-3} and \eqref{Relation-v_44i-3}, equations \eqref{Equality-1st order derivative of F-3} and \eqref{Equality-2nd order derivative of F-3} can be simplified as
\begin{equation}\label{Relation-1st condition-3}
0\sim\sum_{i, j\in G}F^{ij}v_{ij4}
\end{equation}
and
\begin{equation}\label{Relation-2nd condition-3}
-\sum_{i, j}F^{ij}v_{ij44}\sim\sum_{i, j, k, l\in G}F^{ij, kl}v_{ij4}v_{kl4}.
\end{equation}
Substituting \eqref{Relation-2nd condition-3} into \eqref{Equality-F1-3}, we get
\begin{equation}\label{Equality-F2-3}
-\frac{(\det(A))^2}{\sigma_3(G)}\sum_{i, j}F^{ij}\phi_{ij}\sim\sum_{i, j, k, l\in G}F^{ij, kl}v_{ij4}v_{kl4}+2\sum_{p\in G}\sum_{i, j}\frac{1}{v_{pp}}F^{ij}v_{4pi}v_{4pj}.
\end{equation}
 
Next, we carefully analyze each term on the right-hand side of \eqref{Equality-F2-3}. Recall that
\begin{align*}
F(\nabla^2v, \nabla v)
=&~(v_{11}+v_{33})(v_{22}+v_{44})-(v_{12}+v_{34})(v_{21}+v_{43})-(v_{14}-v_{32})(v_{41}-v_{23})\\
&-(v_2^2+v_4^2)(v_{11}+v_{33})-(v_1^2+v_3^2)(v_{22}+v_{44})\\
&+(v_1v_2+v_3v_4)(v_{12}+v_{34}+v_{21}+v_{43})\\
&+(v_1v_4-v_2v_3)(v_{14}-v_{32}+v_{41}-v_{23})-\lambda.
\end{align*}
Direct calculation implies
\begin{align*}
F^{11, 22}&=F^{22, 11}=F^{22, 33}=F^{33, 22}=1, \\
F^{12, 21}&=F^{21, 12}=F^{23, 32}=F^{32, 23}=-1.
\end{align*}
By collecting the terms in \eqref{Equality-F2-3}, we obtain
\begin{equation}\label{Equality-F3-3}
\begin{aligned}
-\frac{(\det(A))^2}{\sigma_3(G)}\sum_{i, j}F^{ij}\phi_{ij} \sim&~\frac{2F^{11}}{v_{11}}v_{114}^2+\frac{2F^{22}}{v_{22}}v_{224}^2+\frac{2F^{33}}{v_{33}}v_{334}^2+\left(\frac{2F^{33}}{v_{11}}+\frac{2F^{11}}{v_{33}}\right)v_{134}^2\\
&+\left(-2+\frac{2F^{22}}{v_{11}}+\frac{2F^{11}}{v_{22}}\right)v_{124}^2+\left(-2+\frac{2F^{33}}{v_{22}}+\frac{2F^{22}}{v_{33}}\right)v_{234}^2\\
&+2v_{114}v_{224}+\frac{4F^{12}}{v_{11}}v_{114}v_{124}+2v_{224}v_{334}+\frac{4F^{12}}{v_{22}}v_{224}v_{124}\\
&+\frac{4F^{23}}{v_{22}}v_{224}v_{234}+\frac{4F^{23}}{v_{33}}v_{334}v_{234}+\frac{4F^{23}}{v_{11}}v_{134}v_{124}+\frac{4F^{12}}{v_{33}}v_{134}v_{234}.
\end{aligned}
\end{equation}
Note that $F^{22}>0$. Since $F^{13}=0$, by \eqref{Relation-1st condition-3}, we have
\begin{equation}\label{Relation-v_224-3}
v_{224}\sim-\frac{F^{11}}{F^{22}}v_{114}-\frac{F^{33}}{F^{22}}v_{334}-\frac{2F^{12}}{F^{22}}v_{124}-\frac{2F^{23}}{F^{22}}v_{234}.
\end{equation}
Putting \eqref{Relation-v_224-3} into \eqref{Equality-F3-3}, we obtain a quadratic form in terms of $v_{114}$, $v_{334}$, $v_{134}$, $v_{124}$ and $v_{234}$, given by
\begin{align*}
-\frac{(\det(A))^2}{2\sigma_3(G)}\sum_{i,j}F^{ij}\phi_{ij}
\sim&\left(\frac{F^{11}}{v_{11}}+\frac{(F^{11})^2}{F^{22}v_{22}}-\frac{F^{11}}{F^{22}}\right)v_{114}^2+\left(\frac{F^{33}}{v_{33}}+\frac{(F^{33})^2}{F^{22}v_{22}}-\frac{F^{33}}{F^{22}}\right)v_{334}^2\\
&+\left(\frac{F^{33}}{v_{11}}+\frac{F^{11}}{v_{33}}\right)v_{134}^2+\left(-1+\frac{F^{22}}{v_{11}}+\frac{F^{11}}{v_{22}}\right)v_{124}^2\\
&+\left(-1+\frac{F^{33}}{v_{22}}+\frac{F^{22}}{v_{33}}\right)v_{234}^2+2\left(\frac{F^{11}F^{33}}{F^{22}v_{22}}-\frac{F^{11}}{2F^{22}}-\frac{F^{33}}{2F^{22}}\right)v_{114}v_{334}\\
&+2\left(\frac{F^{12}}{v_{11}}+\frac{F^{11}F^{12}}{F^{22}v_{22}}-\frac{F^{12}}{F^{22}}\right)v_{114}v_{124}+2\left(\frac{F^{11}F^{23}}{F^{22}v_{22}}-\frac{F^{23}}{F^{22}}\right)v_{114}v_{234}\\
&+2\left(\frac{F^{12}F^{33}}{F^{22}v_{22}}-\frac{F^{12}}{F^{22}}\right)v_{334}v_{124}+2\left(\frac{F^{23}}{v_{33}}+\frac{F^{23}F^{33}}{F^{22}v_{22}}-\frac{F^{23}}{F^{22}}\right)v_{334}v_{234}\\
&+2\frac{F^{23}}{v_{11}}v_{134}v_{124}+2\frac{F^{12}}{v_{33}}v_{134}v_{234}.
\end{align*}
 
To prove inequality \eqref{Inequ-Theorem-minimal rank=3}, we must verify that the above quadratic form in the variables $v_{114}$, $v_{334}$, $v_{134}$, $v_{124}$, and $v_{234}$ is semi-positive definite modulo $\phi$. This requires carefully computing the coefficients of each term and confirming that the corresponding coefficient matrix is semi-positive definite modulo $\phi$.
 
Before arranging these coefficients, we recall the equation $F(\nabla^2v, \nabla v)=0$. It is equivalent to
\begin{equation}\label{Relation-Equation-3}
0\sim (v_{11}+v_{33})v_{22}-(v_2^2+v_4^2)(v_{11}+v_{33})-(v_1^2+v_3^2)v_{22} -\lambda.
\end{equation}
Obviously, the matrix $(F^{ij})$ is given by
\begin{align}\label{Relation-Fij-a}
(F^{ij})=\begin{pmatrix}
F^{11} & F^{12} & 0 & -F^{23}\\[4pt]
F^{12} & F^{22} & F^{23} & 0\\[4pt]
0 & F^{23} & F^{11} & F^{12}\\[4pt]
-F^{23} & 0 & F^{12} & F^{22}
\end{pmatrix},
\end{align}
where
\begin{equation}\label{Relation-Fij-b}
\begin{aligned}
&F^{11}\sim v_{22}-(v_2^2+v_4^2), \hspace{1.65cm} F^{12}=v_1v_2+v_3v_4\\
&F^{22}=(v_{11}+v_{33})-(v_1^2+v_3^2), \quad F^{23}=-(v_1v_4-v_2v_3).
\end{aligned}
\end{equation}
 
We now begin to simplify the coefficients of all terms in the final quadratic form. In the following calculations, we will frequently use \eqref{Relation-Equation-3}--\eqref{Relation-Fij-b}.

$\bullet$ The coefficient of $v_{114}^2$ is
\begin{align*}
\frac{F^{11}}{v_{11}}+\frac{(F^{11})^2}{F^{22}v_{22}}-\frac{F^{11}}{F^{22}}
=&~\frac{1}{F^{22}v_{22}}\frac{F^{11}}{v_{11}}(F^{22}v_{22}+F^{11}v_{11}-v_{11}v_{22})\\
\sim&~\frac{1}{F^{22}v_{22}}\frac{F^{11}}{v_{11}}\left((v_2^2+v_4^2)v_{33}+\lambda\right).
\end{align*}
 
$\bullet$ The coefficient of $v_{134}^2$ is
\begin{align*}
\frac{F^{33}}{v_{11}}+\frac{F^{11}}{v_{33}}
=&~\frac{1}{v_{11}v_{33}}(F^{33}v_{33}+F^{11}v_{11})\\
\sim&~\frac{1}{v_{11}v_{33}}\left((v_1^2+v_3^2)v_{22}+\lambda\right).
\end{align*}
 
$\bullet$ The coefficient of $v_{124}^2$ is
\begin{align*}
-1+\frac{F^{22}}{v_{11}}+\frac{F^{11}}{v_{22}}
=&~\frac{1}{v_{11}v_{22}}(-v_{11}v_{22}+F^{22}v_{22}+F^{11}v_{11})\\
\sim&~\frac{1}{v_{11}v_{22}}\left((v_2^2+v_4^2)v_{33}+\lambda\right).
\end{align*}
 
$\bullet$ The coefficient of $2v_{114}v_{334}$ is
\begin{align*}
\frac{F^{11}F^{33}}{F^{22}v_{22}}-\frac{F^{11}}{2F^{22}}-\frac{F^{33}}{2F^{22}}=&~\frac{1}{F^{22}v_{22}}F^{11}(F^{11}-v_{22})\\
\sim&-\frac{1}{F^{22}v_{22}}F^{11}(v_2^2+v_4^2).
\end{align*}
 
$\bullet$ The coefficient of $2v_{114}v_{124}$ is
\begin{align*}
\frac{F^{12}}{v_{11}}+\frac{F^{11}F^{12}}{F^{22}v_{22}}-\frac{F^{12}}{F^{22}}
=&~\frac{1}{F^{22}v_{22}}\frac{F^{12}}{v_{11}}(F^{22}v_{22}+F^{11}v_{11}-v_{11}v_{22})\\
\sim&~\frac{1}{F^{22}v_{22}}\frac{F^{12}}{v_{11}}\left((v_2^2+v_4^2)v_{33}+\lambda\right).
\end{align*}

$\bullet$ The coefficient of $2v_{114}v_{234}$ is
\begin{align*}
\frac{F^{11}F^{23}}{F^{22}v_{22}}-\frac{F^{23}}{F^{22}}
=&~\frac{1}{F^{22}v_{22}}F^{23}(F^{11}-v_{22})\\
\sim&-\frac{1}{F^{22}v_{22}}F^{23}(v_2^2+v_4^2).
\end{align*}
 
For the remaining coefficients, by exploiting the symmetry between indices $1$ and $3$, we directly obtain the following simplifications.
 
$\bullet$ The coefficient of $v_{334}^2$ is
\begin{align*} \frac{F^{33}}{v_{33}}+\frac{(F^{33})^2}{F^{22}v_{22}}-\frac{F^{33}}{F^{22}}\sim\frac{1}{F^{22}v_{22}}\frac{F^{33}}{v_{33}}\left((v_2^2+v_4^2)v_{11}+\lambda\right).
\end{align*}
 
$\bullet$ The coefficient of $v_{234}^2$ is
\begin{align*}
-1+\frac{F^{33}}{v_{22}}+\frac{F^{22}}{v_{33}}\sim\frac{1}{v_{22}v_{33}}\left((v_2^2+v_4^2)v_{11}+\lambda\right).
\end{align*}

$\bullet$ The coefficient of $2v_{334}v_{124}$ is
\begin{align*}
\frac{F^{12}F^{33}}{F^{22}v_{22}}-\frac{F^{12}}{F^{22}}\sim-\frac{1}{F^{22}v_{22}}F^{12}(v_2^2+v_4^2).
\end{align*}
 
$\bullet$ The coefficient of $2v_{334}v_{234}$ is
\begin{align*}
\frac{F^{23}}{v_{33}}+\frac{F^{23}F^{33}}{F^{22}v_{22}}-\frac{F^{23}}{F^{22}}\sim\frac{1}{F^{22}v_{22}}\frac{F^{23}}{v_{33}}\left((v_2^2+v_4^2)v_{11}+\lambda\right).
\end{align*}
 
$\bullet$ The coefficients of $2v_{114}v_{134}$, $2v_{334}v_{134}$, $2v_{134}v_{124}$, $2v_{134}v_{234}$ and $2v_{124}v_{234}$ are
\begin{align*}
0, \quad 0, \quad \frac{F^{23}}{v_{11}}, \quad \frac{F^{12}}{v_{33}} \ \mathrm{and} \ 0,
\end{align*}
respectively.
 
Therefore, we obtain 
\begin{align}\label{Quadratic form-f}
-\frac{(\det(A))^2}{2\sigma_3(G)} \sum_{i, j}F^{ij}\phi_{ij}\sim f(X),
\end{align}
where $ f(X)=X^TMX$ is a quadratic form with $X=(v_{114}, v_{334}, v_{134}, v_{124}, v_{234})^T$ as the associated vector and $M$ as the symmetric matrix defined below.  
\begin{center}
$\begin{matrix}
 & v_{114} & v_{334} & v_{134} & v_{124} & v_{234}\\[10pt]
v_{114} & \makecell{\frac{1}{F^{22}v_{22}}\frac{F^{11}}{v_{11}} \cdot((v_2^2\\[4pt] +v_4^2)v_{33}+\lambda)} & \makecell{-\frac{1}{F^{22}v_{22}}F^{11}\\[4pt] \cdot (v_2^2+v_4^2)} & 0 & \makecell{\frac{1}{F^{22}v_{22}}\frac{F^{12}}{v_{11}}\cdot((v_2^2\\[4pt] +v_4^2)v_{33}+\lambda)} & \makecell{-\frac{1}{F^{22}v_{22}}F^{23}\\[4pt] \cdot (v_2^2+v_4^2)}\\[18pt]
v_{334} & \makecell{-\frac{1}{F^{22}v_{22}}F^{11}\\[4pt] \cdot (v_2^2+v_4^2)} & \makecell{\frac{1}{F^{22}v_{22}}\frac{F^{33}}{v_{33}}\cdot((v_2^2\\[4pt] +v_4^2)v_{11}+\lambda)} & 0 & \makecell{-\frac{1}{F^{22}v_{22}}F^{12}\\[4pt] \cdot (v_2^2+v_4^2)} & \makecell{\frac{1}{F^{22}v_{22}}\frac{F^{23}}{v_{33}}\cdot((v_2^2\\[4pt] +v_4^2)v_{11}+\lambda)}\\[18pt]
v_{134} & 0 & 0 & \makecell{\frac{1}{v_{11}v_{33}}\cdot((v_1^2\\[4pt] +v_3^2)v_{22}+\lambda)} & \frac{F^{23}}{v_{11}} & \frac{F^{12}}{v_{33}}\\[18pt]
v_{124} & \makecell{\frac{1}{F^{22}v_{22}}\frac{F^{12}}{v_{11}}\cdot((v_2^2\\[4pt] +v_4^2)v_{33}+\lambda)} & \makecell{-\frac{1}{F^{22}v_{22}}F^{12}\\[4pt] \cdot (v_2^2+v_4^2)} & \frac{F^{23}}{v_{11}} & \makecell{\frac{1}{v_{11}v_{22}}\cdot((v_2^2\\[4pt] +v_4^2)v_{33}+\lambda)} & 0\\[18pt]
v_{234} & \makecell{-\frac{1}{F^{22}v_{22}}F^{23}\\[4pt] \cdot (v_2^2+v_4^2)} & \makecell{\frac{1}{F^{22}v_{22}}\frac{F^{23}}{v_{33}}\cdot((v_2^2\\[4pt] +v_4^2)v_{11}+\lambda)} & \frac{F^{12}}{v_{33}} & 0 & \makecell{\frac{1}{v_{22}v_{33}}\cdot((v_2^2\\[4pt] +v_4^2)v_{11}+\lambda)}
\end{matrix}$
\end{center}
 
\vspace{0.4cm} Let us verify that the matrix $M$ is semi-positive modulo $\phi$, i.e., the corresponding quadratic form $f(X)$ satisfies $f(X)\gtrsim0$. This is highly technical and challenging. 
 
First, we perform congruence transformations on the matrix $M$, which preserve the positivity of $M$, to convert it into a block diagonal form. The specific transformations are as follows.
 
\begin{enumerate}
\item The first column is multiplied by $-\frac{F^{12}}{F^{11}}$ and added to the fourth column, while the first row is multiplied by $-\frac{F^{12}}{F^{11}}$ and added to the fourth row.
  
\item The second column is multiplied by $-\frac{F^{23}}{F^{11}}$ and added to the fifth column, while the second row is multiplied by $-\frac{F^{23}}{F^{11}}$ and added to the fifth row.
\end{enumerate}
It yields the following block diagonal matrix.
\begin{center}
$\begin{matrix}
 & v_{114} & v_{334} & v_{134} & v_{124} & v_{234}\\[10pt]
v_{114} & \makecell{\frac{1}{F^{22}v_{22}}\frac{F^{11}}{v_{11}} \cdot((v_2^2\\[4pt] +v_4^2)v_{33}+\lambda)} & \makecell{-\frac{1}{F^{22}v_{22}}F^{11}\\ \cdot (v_2^2+v_4^2)} & 0 & 0 & 0\\[18pt]
v_{334} & \makecell{-\frac{1}{F^{22}v_{22}}F^{11}\\ \cdot (v_2^2+v_4^2)} & \makecell{\frac{1}{F^{22}v_{22}}\frac{F^{33}}{v_{33}}\cdot((v_2^2\\[4pt] +v_4^2)v_{11}+\lambda)} & 0 & 0 & 0\\[18pt]
v_{134} & 0 & 0 & \makecell{\frac{1}{v_{11}v_{33}}\cdot((v_1^2\\[4pt] +v_3^2)v_{22}+\lambda)} & \frac{F^{23}}{v_{11}} & \frac{F^{12}}{v_{33}}\\[18pt]
v_{124} & 0 & 0 & \frac{F^{23}}{v_{11}} & \makecell{\frac{F^{11}F^{22}-(F^{12})^2}{F^{11}v_{11}\cdot F^{22}v_{22}}\cdot((v_2^2\\[4pt] +v_4^2)v_{33}+\lambda)} & \makecell{\frac{1}{F^{22}v_{22}}\frac{F^{12}F^{23}}{F^{11}}\\[4pt] \cdot(v_2^2+v_4^2)}\\[18pt]
v_{234} & 0 & 0 & \frac{F^{12}}{v_{33}} & \makecell[l]{\frac{1}{F^{22}v_{22}}\frac{F^{12}F^{23}}{F^{11}}\\[4pt] \cdot(v_2^2+v_4^2)} & \makecell{\frac{F^{22}F^{33}-(F^{23})^2}{F^{22}v_{22}\cdot F^{33}v_{33}}((v_2^2\\[4pt] +v_4^2)v_{11}+\lambda)}
\end{matrix}$
\end{center}

\textbf{Claim 2. } The matrix $$M_1=
\begin{pmatrix}
\makecell{\frac{1}{F^{22}v_{22}}\frac{F^{11}}{v_{11}} \cdot((v_2^2\\[4pt] +v_4^2)v_{33}+\lambda)} & \makecell{-\frac{1}{F^{22}v_{22}}F^{11}\\ \cdot (v_2^2+v_4^2)}\\[18pt]
\makecell{-\frac{1}{F^{22}v_{22}}F^{11}\\ \cdot (v_2^2+v_4^2)} & \makecell{\frac{1}{F^{22}v_{22}}\frac{F^{33}}{v_{33}}\cdot((v_2^2\\[4pt] +v_4^2)v_{11}+\lambda)}
\end{pmatrix}$$ is positive definite modulo $\phi$.
 
\textbf{Claim 3. } The matrix $$M_2=
\begin{pmatrix}
\makecell{\frac{1}{v_{11}v_{33}}\cdot((v_1^2\\[4pt] +v_3^2)v_{22}+\lambda)} & \frac{F^{23}}{v_{11}} & \frac{F^{12}}{v_{33}}\\[18pt]
\frac{F^{23}}{v_{11}} & \makecell{\frac{F^{11}F^{22}-(F^{12})^2}{F^{11}v_{11}\cdot F^{22}v_{22}}\cdot((v_2^2\\[4pt] +v_4^2)v_{33}+\lambda)} & \makecell{\frac{1}{F^{22}v_{22}}\frac{F^{12}F^{23}}{F^{11}}\\[4pt] \cdot(v_2^2+v_4^2)}\\[18pt]
\frac{F^{12}}{v_{33}} & \makecell{\frac{1}{F^{22}v_{22}}\frac{F^{12}F^{23}}{F^{11}}\\[4pt] \cdot(v_2^2+v_4^2)} & \makecell{\frac{F^{22}F^{33}-(F^{23})^2}{F^{22}v_{22}\cdot F^{33}v_{33}}((v_2^2\\[4pt] +v_4^2)v_{11}+\lambda)}
\end{pmatrix}$$ is positive definite modulo $\phi$.
 
Once \textbf{Claim 2} and \textbf{Claim 3} are proved, we will know that the symmetric matrix $M$ is positive definite modulo $\phi$. That is, the corresponding quadratic form $f(X)$ satisfies $$f(X)=X^TMX\gtrsim0. $$ Combining this with formula \eqref{Quadratic form-f}, we conclude the proof of Theorem \ref{Theorem-minimal rank is 3}.
\end{proof}

Our remaining task is to prove \textbf{Claim 2} and \textbf{Claim 3}.

\textbf{Proof of Claim 2.} Let us compute the leading principal minors of $M_1$, denoted by $P_1(M_1)$ and  $P_2(M_1)$.

$\bullet$ $P_1(M_1)$. It is clear that $$P_1(M_1)=\frac{1}{F^{22}v_{22}} \frac{F^{11}}{v_{11}}\left((v_2^2+v_4^2)v_{33}+\lambda\right)>0. $$
 
$\bullet$ $P_2(M_1)$. Using \eqref{Relation-Equation-3} and \eqref{Relation-Fij-b}, we obtain 
\begin{align}\label{Relation-F22v22}
F^{22}v_{22}\sim(v_2^2+v_4^2)(v_{11}+v_{33})+\lambda. 
\end{align} 
Note that $F^{33}=F^{11}$. Then we have
\begin{align*}
P_2(M_1)=&~\frac{1}{F^{22}v_{22}}\frac{F^{11}}{v_{11}}\left((v_2^2+v_4^2)v_{33}+\lambda\right)\cdot\frac{1}{F^{22}v_{22}}\frac{F^{33}}{v_{33}}\left((v_2^2+v_4^2)v_{11}+\lambda\right)\\
&-\frac{1}{(F^{22}v_{22})^2}(F^{11})^2(v_2^2+v_4^2)^2\\
=&~\frac{\lambda}{(F^{22}v_{22})^2}\frac{(F^{11})^2}{v_{11}v_{33}}\left((v_2^2+v_4^2)(v_{11}+v_{33})+\lambda\right)\\
\sim&~\frac{\lambda}{F^{22}v_{22}}\frac{(F^{11})^2}{v_{11}v_{33}}>0.
\end{align*}
Therefore, the matrix $M_1$ is positive definite modulo $\phi$. \qed

\vspace{0.4cm}\textbf{Proof of Claim 3.} Along the same route, we compute the leading principal minors of $M_2$, denoted by $P_1(M_2)$,  $P_2(M_2)$ and  $P_3(M_2)$.

$\bullet$ $P_1(M_2)$. Obviously, $$P_1(M_2)=\frac{1}{v_{11}v_{33}} \left((v_1^2+v_3^2)v_{22}+\lambda\right)>0. $$

$\bullet$ $P_2(M_2)$. Next, we proceed to compute $P_2(M_2)$ as follows
\begin{align*}
P_2(M_2)=&~\frac{1}{v_{11}v_{33}}\left((v_1^2+v_3^2)v_{22}+\lambda\right)\cdot\frac{F^{11}F^{22}-(F^{12})^2}{F^{11}v_{11}\cdot F^{22}v_{22}}\left((v_2^2+v_4^2)v_{33}+\lambda\right)-\left(\frac{F^{23}}{v_{11}}\right)^2\\
=&~\frac{1}{F^{11}v_{11}\cdot F^{22}v_{22}}\frac{1}{v_{11}v_{33}}\\
&\cdot\bigg(\left(F^{11}F^{22}-(F^{12})^2\right)\left((v_1^2+v_3^2)v_{22}+\lambda\right)\left((v_2^2+v_4^2)v_{33}+\lambda\right)-F^{22}v_{22}\cdot F^{33}v_{33}(F^{23})^2\bigg).
\end{align*}
From \eqref{Relation-Equation-3} and the expressions for $F^{11}$, $F^{22}$ and $F^{33}$ in \eqref{Relation-Fij-a}--\eqref{Relation-Fij-b}, we directly have
\begin{align*}
&\left((v_1^2+v_3^2)v_{22}+\lambda\right)\left((v_2^2+v_4^2)v_{33}+\lambda\right)\\
=&~(v_1^2+v_3^2)(v_2^2+v_4^2)v_{22}v_{33}+\lambda\left((v_1^2+v_3^2)v_{22}+(v_2^2+v_4^2)v_{33}+\lambda\right)\\
\sim&~(v_1^2+v_3^2)(v_2^2+v_4^2)v_{22}v_{33}+\lambda\left((v_{11}+v_{33})v_{22}-(v_2^2+v_4^2)v_{11}\right)\\
\sim&~(v_1^2+v_3^2)(v_2^2+v_4^2)v_{22}v_{33}+\lambda\cdot v_{22}v_{33}+\lambda\cdot F^{11}v_{11}\\
\sim&~F^{22}v_{22}\cdot F^{33}v_{33}+\lambda\cdot F^{11}v_{11}.
\end{align*}
It yields
\begin{align*}
P_2(M_2)\sim&~\frac{1}{F^{11}v_{11}\cdot F^{22}v_{22}}\frac{1}{v_{11}v_{33}}\cdot\bigg(F^{22}v_{22}\cdot F^{33}v_{33}\left(F^{11}F^{22}-(F^{12})^2\right)\\
&\hspace{0.5cm}+\lambda\cdot F^{11}v_{11}\left(F^{11}F^{22} -(F^{12})^2\right)-F^{22}v_{22}\cdot F^{33}v_{33}(F^{23})^2\bigg)\\
=&~\frac{1}{F^{11}v_{11}\cdot F^{22}v_{22}}\frac{1}{v_{11}v_{33}}F^{22}v_{22}\cdot F^{33}v_{33}\left(F^{11}F^{22}-(F^{12})^2-(F^{23})^2\right)\\
&+\frac{\lambda}{F^{22}v_{22}}\frac{1}{v_{11}v_{33}}\left(F^{11}F^{22}-(F^{12})^2\right).
\end{align*}
Once again, by \eqref{Relation-Equation-3} and \eqref{Relation-Fij-b}, we have
\begin{equation}\label{Relation-elliptic coefficients-3}
\begin{aligned}
F^{11}F^{22}-(F^{12})^2-(F^{23})^2
\sim&~(v_{11}+v_{33})v_{22}-(v_2^2+v_4^2)(v_{11}+v_{33})-(v_1^2+v_3^2)v_{22}\\
&+(v_1^2+v_3^2)(v_2^2+v_4^2)\\
&-(v_1v_2+v_3v_4)^2-(v_1v_4-v_2v_3)^2\\
=&~(v_{11}+v_{33})v_{22}-(v_2^2+v_4^2)(v_{11}+v_{33})-(v_1^2+v_3^2)v_{22}\\
\sim&~\lambda.
\end{aligned}
\end{equation}
So we obtain
\begin{align*}
P_2(M_2)\sim&~\frac{\lambda}{F^{11}v_{11}\cdot F^{22}v_{22}}\frac{1}{v_{11}v_{33}}F^{22}v_{22}\cdot F^{33}v_{33}+\frac{\lambda}{F^{22}v_{22}}\frac{1}{v_{11}v_{33}}\left(\lambda+(F^{23})^2\right)>0.
\end{align*}
 
$\bullet$ $P_3(M_2)$. By \eqref{Relation-F22v22}, we know that $$((v_2^2+v_4^2)v_{33}+\lambda)((v_2^2+v_4^2)v_{11}+\lambda)\sim (v_2^2+v_4^2)^2v_{11}v_{33}+F^{22}v_{22}. $$ Noting that $F^{33}=F^{11}$, we  can derive
\begin{align*}
P_3(M_2)=&~\frac{1}{v_{11}v_{33}}\left((v_1^2+v_3^2)v_{22}+\lambda\right)\\
&\cdot\bigg(\frac{(F^{11}F^{22}-(F^{12})^2)(F^{22}F^{33}-(F^{23})^2)}{F^{11}v_{11}\cdot(F^{22}v_{22})^2\cdot F^{33}v_{33}}\left((v_2^2+v_4^2)v_{33}+\lambda\right)\left((v_2^2+v_4^2)v_{11}+\lambda\right)\\
&\hspace{0.6cm}-\frac{1}{(F^{22}v_{22})^2}\frac{(F^{12})^2(F^{23})^2}{(F^{11})^2}(v_2^2+v_4^2)^2\bigg)\\
&~-\frac{F^{23}}{v_{11}}\bigg(\frac{F^{23}}{v_{11}}\frac{F^{22}F^{33}-(F^{23})^2}{F^{22}v_{22}\cdot F^{33}v_{33}}\left((v_2^2+v_4^2)v_{11}+\lambda\right)-\frac{F^{12}}{v_{33}}\frac{1}{F^{22}v_{22}}\frac{F^{12}F^{23}}{F^{11}}(v_2^2+v_4^2)\bigg)\\
&~+\frac{F^{12}}{v_{33}}\bigg(\frac{F^{23}}{v_{11}}\frac{1}{F^{22}v_{22}}\frac{F^{12}F^{23}}{F^{11}}(v_2^2+v_4^2)-\frac{F^{12}}{v_{33}}\frac{F^{11}F^{22}-(F^{12})^2}{F^{11}v_{11}\cdot F^{22}v_{22}}\left((v_2^2+v_4^2)v_{33}+\lambda\right)\bigg)\\
=&~\frac{1}{v_{11}v_{33}}\left((v_1^2+v_3^2)v_{22}+\lambda\right)\\
&\cdot\bigg(\frac{F^{11}F^{22}(F^{11}F^{22}-(F^{12})^2-(F^{23})^2)}{F^{11}v_{11}\cdot(F^{22}v_{22})^2\cdot  F^{33}v_{33}}\left((v_2^2+v_4^2)v_{33}+\lambda\right)\left((v_2^2+v_4^2)v_{11}+\lambda\right)\\
&\hspace{0.6cm}+\frac{(F^{12})^2(F^{23})^2}{F^{11}v_{11}\cdot(F^{22}v_{22})^2\cdot F^{33}v_{33}}\left((v_2^2+v_4^2)v_{11}+\lambda\right)\left((v_2^2+v_4^2)v_{33}+\lambda\right)\\
&\hspace{0.6cm}-\frac{1}{(F^{22}v_{22})^2}\frac{(F^{12})^2(F^{23})^2}{(F^{11})^2}(v_2^2+v_4^2)^2\bigg)\\
&~-\frac{1}{v_{11}}\frac{(F^{23})^2}{F^{22}v_{22}\cdot F^{33}v_{33}}(v_2^2+v_4^2)\bigg(F^{11}F^{22}-(F^{12})^2-(F^{23})^2\bigg)\\
&-\frac{\lambda}{v_{11}^2}\frac{(F^{23})^2}{F^{22}v_{22}\cdot F^{33}v_{33}}\hspace{1pt}\bigg(F^{22}F^{33}-(F^{23})^2\bigg)\\
&~-\frac{1}{v_{33}}\frac{(F^{12})^2}{F^{11}v_{11}\cdot F^{22}v_{22}}(v_2^2+v_4^2)\bigg(F^{11}F^{22}-(F^{12})^2-(F^{23})^2\bigg)\\
&~-\frac{\lambda}{v_{33}^2}\frac{(F^{12})^2}{F^{11}v_{11}\cdot F^{22}v_{22}}\bigg(F^{11}F^{22}-(F^{12})^2\bigg)\\
\sim&~\frac{F^{11}F^{22}}{F^{11}v_{11}\cdot (F^{22}v_{22})^2 \cdot F^{33}v_{33}}(v_2^2+v_4^2)^2\left((v_1^2+v_3^2)v_{22}+\lambda\right)\bigg(F^{11}F^{22}-(F^{12})^2-(F^{23})^2\bigg)\\
&~+\frac{\lambda}{v_{11}v_{33}}\frac{F^{11}F^{22}}{F^{11}v_{11}\cdot F^{22}v_{22}\cdot F^{33}v_{33}}\left((v_1^2+v_3^2)v_{22}+\lambda\right)\bigg(F^{11}F^{22}-(F^{12})^2-(F^{23})^2\bigg)\\
&~+\frac{\lambda}{v_{11}v_{33}}\frac{(F^{12})^2(F^{23})^2}{F^{11}v_{11}\cdot F^{22}v_{22}\cdot F^{33}v_{33}}\left((v_1^2+v_3^2)v_{22}+\lambda\right)\\
&~-\frac{\lambda}{v_{33}}\frac{1}{F^{11}v_{11}\cdot F^{22}v_{22}}\bigg(\frac{(F^{23})^2}{v_{11}}\left(F^{22}F^{33}-(F^{23})^2\right)+\frac{(F^{12})^2}{v_{33}}\left(F^{11}F^{22}-(F^{12})^2\right)\bigg)\\
&~-\frac{1}{v_{33}}\frac{(F^{12})^2+(F^{23})^2}{F^{11}v_{11}\cdot F^{22}v_{22}}(v_2^2+v_4^2)\bigg(F^{11}F^{22}-(F^{12})^2-(F^{23})^2\bigg)\\
\stackrel{\triangle}{=}&~\circled{1}+\circled{2}+\circled{3}+\circled{4}+\circled{5}.
\end{align*}
 
We now group the five terms above into two parts and evaluate each group separately. Since $$(F^{12})^2+(F^{23})^2=(v_1^2+v_3^2)(v_2^2+v_4^2), $$ we have
\begin{equation}\label{Equality-2+5}
\circled{1}+\circled{5}=\frac{\lambda\hspace{1pt}F^{11}F^{22}}{F^{11}v_{11}\cdot (F^{22}v_{22})^2\cdot F^{33}v_{33}}(v_2^2+v_4^2)^2\bigg(F^{11}F^{22}-(F^{12})^2-(F^{23})^2\bigg).
\end{equation}
It is evident that
\begin{align*}
&~\circled{2}+\circled{3}+\circled{4}\\
=&~\frac{\lambda}{v_{11}v_{33}}\frac{1}{F^{11}v_{11}\cdot F^{22}v_{22}\cdot F^{33}v_{33}}\\
&\cdot\bigg(\left((v_1^2+v_3^2)v_{22}+\lambda\right)F^{11}F^{22}\left(F^{11}F^{22}-(F^{12})^2-(F^{23})^2\right)\\
&\hspace{0.6cm}+\left((v_1^2+v_3^2)v_{22}+\lambda\right)(F^{12})^2(F^{23})^2\\
&\hspace{0.6cm}-F^{33}v_{33}(F^{23})^2\left(F^{22}F^{33}-(F^{23})^2\right)-F^{11}v_{11}(F^{12})^2\left(F^{11}F^{22}-(F^{12})^2\right)\bigg)
\end{align*}
From relations \eqref{Relation-elliptic coefficients-3} and \eqref{Relation-Equation-3}--\eqref{Relation-Fij-b}, we obtain
\begin{align*}
&(F^{23})^2\leq F^{11}F^{22}-(F^{12})^2, \\
&(F^{12})^2\leq F^{22}F^{33}-(F^{23})^2, \\
\intertext{and}
&F^{11}v_{11}+F^{33}v_{33}\sim (v_1^2+v_3^2)v_{22}+\lambda,
\end{align*}
which lead to
\begin{align*}
&-F^{33}v_{33}(F^{23})^2\left(F^{22}F^{33}-(F^{23})^2\right)-F^{11}v_{11}(F^{12})^2\left(F^{11}F^{22}-(F^{12})^2\right)\\
\geq&-\left(F^{11}v_{11}+F^{33}v_{33}\right)\left(F^{11}F^{22}-(F^{12})^2\right)\left(F^{22}F^{33}-(F^{23})^2\right)\\
\sim&-\left((v_1^2+v_3^2)v_{22}+\lambda\right)\left(F^{11}F^{22}-(F^{12})^2\right)\left(F^{22}F^{33}-(F^{23})^2\right).
\end{align*}
Therefore, 
\begin{equation}\label{Inequality-1+3+4}
\begin{aligned}
\circled{2}+\circled{3}+\circled{4}
\gtrsim&~\frac{\lambda}{v_{11}v_{33}}\frac{1}{F^{11}v_{11}\cdot F^{22}v_{22}\cdot F^{33}v_{33}}\left((v_1^2+v_3^2)v_{22}+\lambda\right)\\
&\cdot\bigg(F^{11}F^{22}\left(F^{11}F^{22}-(F^{12})^2-(F^{23})^2\right)+(F^{12})^2(F^{23})^2\\
&\hspace{0.6cm}-\left(F^{11}F^{22}-(F^{12})^2\right)\left(F^{22}F^{33}-(F^{23})^2\right)\bigg)\\
=&~0.
\end{aligned}
\end{equation}
Combining \eqref{Equality-2+5}--\eqref{Inequality-1+3+4} and using \eqref{Relation-elliptic coefficients-3}, we conclude that 
\begin{align*}
P_3(M_2)\gtrsim&~\frac{\lambda\hspace{1pt}F^{11}F^{22}}{F^{11}v_{11}\cdot (F^{22}v_{22})^2\cdot F^{33}v_{33}}(v_2^2+v_4^2)^2\bigg(F^{11}F^{22}-(F^{12})^2-(F^{23})^2\bigg)\\
\sim&~\frac{\lambda^2\hspace{1pt}F^{11}F^{22}}{F^{11}v_{11}\cdot (F^{22}v_{22})^2\cdot F^{33}v_{33}}(v_2^2+v_4^2)^2\geq0.
\end{align*}
So the matrix $M_2$ is semi-positive definite modulo $\phi$. \qed

\vspace{5pt} At the end of this section, we will use Theorems \ref{Theorem-minimal rank is 2} and \ref{Theorem-minimal rank is 3} to establish the constant rank theorem.
\begin{corollary}[Constant rank theorem]\label{Corollary}
Let $\Omega\subset\mathbb{R}^4$ be a domain, and let $v\in C^4(\Omega)$ be a convex solution of \eqref{Equation-real equation(v)}. Then the Hessian $\nabla^2v$ has constant rank throughout $\Omega$.
\end{corollary}
\begin{proof}
Define $$l=\min_{x\in\Omega}\mathrm{rank}(\nabla^2v(x)). $$ As discussed in the Section \ref{Section 2}, $l$ can only be $2$ or $3$. Suppose there exists some point $x_0\in\Omega$ such that $\mathrm{rank}(\nabla^2v(x_0))=l$. We will show that $\mathrm{rank}(\nabla^2v(x))=l$ for all $x\in\Omega$. By \eqref{rank identity}, namely $$\mathrm{rank}(\nabla^2v(x))=\mathrm{rank}(K(x))+2,$$ it suffices to study the constant rank property of $K(x)$. Define the set $$\Omega^\ast=\{x\in\Omega: \mathrm{rank}(K(x))=l-2\}. $$ It is clear that $\Omega^\ast$ is relatively closed in $\Omega$. By the strong maximum principle, Theorems \ref{Theorem-minimal rank is 2} and \ref{Theorem-minimal rank is 3} imply that $\Omega^\ast$ is also open. Hence, $\Omega^\ast=\Omega$, meaning that $K$ has constant rank $l-2$ on all of $\Omega$, which in turn implies that $\nabla^2v$ has constant rank $l$ throughout $\Omega$.

This completes the proof of Corollary \ref{Corollary}.
\end{proof}

\section{Proof of Theorem \ref{Theorem-Main}}\label{Section 4}

In this section, we present the proof of Theorem \ref{Theorem-Main} using the continuity method. Based on Corollary \ref{Corollary}, the argument follows a standard approach, as outlined by Caffarelli and Friedman \cite{Caffarelli-Friedman1985}, and further developed in \cite{Korevaar-Lewis1987, Liu-Ma-Xu2010, Ma-Xu2008, Zhang-Zhou2023, Zhang-Zhou2025}.

First, we need to show that the solution of \eqref{Equation-Dirichlet problem-real} is strictly convex when the domain $\Omega$ is a ball. For completeness, we include a proof along the lines of McCuan \cite{McCuan2002}. Similar results can also be found in \cite{Huang2019, Liu-Ma-Xu2010}.
\begin{lemma}\label{Lemma-ball}
Let $B_R(0)$ be the ball in $\mathbb{R}^4$ centered at the origin with radius $R>0$. Let $u\in C^\infty(B_R)\cap C^{1, 1}(\bar{B}_R)$ be the unique solution to the eigenvalue problem \eqref{Equ:EVP-u} in $B_R(0)$. Then $v=-\log(-u/4)$ is strictly convex in $B_R(0)$.
\end{lemma}
\begin{proof}
By uniqueness, the solution $u$ must be radially symmetric when $\Omega$ is a ball. Thus, we can express it as $$u(x)=u(|x|)=u(r), $$ where $r=|x|\in[0,R]$. In this case, $u(r)<0$ for $r\in[0, R)$, and it satisfies the boundary conditions $u'(0)=0$ and $u(R)=0$.
 
A straightforward computation yields $$u_{ij}=\left(r^{-2}u''-r^{-3}u'\right)x_ix_j+r^{-1}u'\delta_{ij}. $$ Notice that equation \eqref{Equation-complexMA(eigenvalue)} can be written in real coordinates as $$16\lambda u^2=(u_{11}+u_{33})(u_{22}+u_{44})-(u_{12}+u_{34})(u_{21}+u_{43})-(u_{14}-u_{32})(u_{41}-u_{23}).$$ It follows that
\begin{equation}\label{Equation-radially symmetric-u}
16\lambda u^2=2r^{-1}u'u''+2r^{-2}(u')^2 \quad \mathrm{in} \ (0, R).
\end{equation}
Since $u$ is plurisubharmonic, we have $\Delta u>0$, which is equivalent to 
\begin{align}\label{Equality-u''}
u''+3r^{-1}u'>0
\end{align} 
for $r\in(0, R)$. In particular, taking the limit as $r\rightarrow0+$ and applying L'H{\^o}pital's rule yields $u''(0)\geq0$. Moreover, by evaluating both sides of equation \eqref{Equation-radially symmetric-u} as $r\rightarrow0+$, we further obtain $u''(0)>0$. Since $v=-\log(-u/4)$, we have 
\begin{align*}
&v'(0)=-\frac{u'(0)}{u(0)}=0, \\
&v''(0)=-\frac{u''(0)}{u(0)}+\left(\frac{u'(0)}{u(0)}\right)^2=-\frac{u''(0)}{u(0)}>0.
\end{align*}
 
We now aim to show that $v''(r)>0$ for all $r\in[0, R)$. Suppose that there exists some $r_0\in(0, R)$ such that $v''(r)>0$ for $r\in[0, r_0)$ but $v''(r_0)=0$. It is clear that $v'''(r_0)\leq0$. On the other hand, $v$ satisfies the following ordinary differential equation 
\begin{equation}\label{Equation-radially symmetric-v}
rv'v''-r(v')^3+(v')^2=8\lambda r^2\quad \mathrm{in} \ (0, R).
\end{equation} 
Differentiating both sides of equation \eqref{Equation-radially symmetric-v} and evaluating at $r_0$, we obtain 
\begin{align}\label{Equality-v'''}
v'''(r_0)=\frac{16\lambda}{v'(r_0)}+\frac{(v'(r_0))^2}{r_0}. 
\end{align} 
By \eqref{Equality-u''}, we have $$v''-(v')^2+3r^{-1}v'>0$$ for $r\in(0, R)$. Since $v''(r_0)=0$, it follow that $v'(r_0)>0$. Combining this with \eqref{Equality-v'''}, we see that $v'''(r_0)>0$, which contradicts the earlier conclusion that $v'''(r_0)\leq0$. Therefore, $v''(r)>0$ for all $r\in[0, R)$, and $v$ is strictly convex in $B_R(0)$.
\end{proof}

Next, we also require the convexity estimates near the boundary for the function $v=-\log(-u/4)$, where $u$ solves the eigenvalue problem \eqref{Equ:EVP-u}.
\begin{lemma}[Convexity estimates near the boundary, see Korevaar \cite{Korevaar1983}, Caffarelli and Friedman \cite{Caffarelli-Friedman1985}]\label{Lemma-BCE}
Let $\Omega\subset\mathbb{R}^n$ be a bounded, smooth, strictly convex domain, and let $u\in C^{1, 1}(\bar{\Omega})\cap C^{\infty}(\Omega)$ satisfy $$u<0\quad \mathrm{in}\ \Omega, \quad u=0\quad \mathrm{on}\ \partial\Omega, \quad |\nabla u|>0\quad \mathrm{on}\ \partial\Omega. $$ Then there exists $\varepsilon>0$ such that the function $v=-\log(-u/4)$ is strictly convex in the boundary strip $\Omega\setminus\bar{\Omega}_\varepsilon$, where
$$\Omega_\varepsilon=\{x\in\Omega: \mathrm{d}(x, \partial\Omega)>\varepsilon\}. $$
\end{lemma}

Finally, we turn to the proof of Theorem \ref{Theorem-Main}.
\begin{proof}[Proof of Theorem \ref{Theorem-Main}]
If $\Omega=B_1$, the unit ball in $\mathbb{R}^4$, then by Lemma \ref{Lemma-ball}, the unique solution to the Dirichlet problem \eqref{Equation-Dirichlet problem-real} is strictly convex. For an arbitrary bounded, smooth, strictly convex domain $\Omega_1=\Omega$, define $$\Omega_t=(1-t)B_1+t\Omega\quad \mathrm{for}\ t\in(0, 1). $$ By the theory of convex bodies, the family $\{\Omega_t\}_{t\in(0, 1)}$ provides a smooth deformation of $B_1$ into $\Omega$, preserving boundedness, smoothness, and strict convexity. For further details, see Schneider's excellent book \cite{Schneider2014}.
 
For $t\in(0,1)$, assume that the function $v_t$ solves the eigenvalue problem
$$\begin{cases}
F(\nabla^2v, \nabla v)=0 & \mathrm{in}\ \Omega_t, \\
\hspace{1.25cm}v(x)\rightarrow+\infty & \mathrm{as}\ \Omega_t\ni x\rightarrow\partial\Omega_t, \\
\hspace{0.9cm}\inf_{\Omega_t}v=2\log2,
\end{cases}$$ where
\begin{align*}
F(\nabla^2v, \nabla v) =&~(v_{11}+v_{33})(v_{22}+v_{44})-(v_{12}+v_{34})(v_{21}+v_{43})-(v_{14}-v_{32})(v_{41}-v_{23})\\
&-(v_2^2+v_4^2)(v_{11}+v_{33})-(v_1^2+v_3^2)(v_{22}+v_{44})\\
&+(v_1v_2+v_3v_4)(v_{12}+v_{34}+v_{21}+v_{43})\\
&+(v_1v_4-v_2v_3)(v_{14}-v_{32}+v_{41}-v_{23})-\lambda.
\end{align*}
According to the a priori estimates for the eigenvalue problem \eqref{Equ:EVP-u} established by Chu, Liu and McCleerey \cite{Chu-Liu-McCleerey-arXiv2024}, we have uniform bounds for $||v_t||_{C^3(\Omega_t)}$ and $||v_t||_{C^{1, 1}(\bar\Omega_t)}$, depending only on the geometry of $\Omega$ and independent of $t$. On the one hand, since $v_0$ is strictly convex, it follows that $v_t$ remains strictly convex for all $t\in(0, \delta)$ for some sufficiently small $\delta>0$. On the other hand, if $v_t$ is strictly convex in $\Omega_t$ for all $t\in(0, t^\ast)$ with $t^\ast\in(0, 1)$, then $v_{t^\ast}$ must be convex in $\Omega_{t^\ast}$.
 
Now, suppose that $v_1$ is not strictly convex in $\Omega_1$. Then there exists some $t_0\in(0, 1)$ at which $v_{t_0}$ is convex but not strictly convex for the first time. By Corollary \ref{Corollary} (Constant rank theorem), the Hessian $\nabla^2v_{t_0}$ is degenerate at every point of $\Omega_{t_0}$. However, Lemma \ref{Lemma-BCE} (Convexity estimates near the boundary) implies that $\nabla^2v_{t_0}$ is of full rank in the boundary strip $\Omega_{t_0}\setminus\bar{\Omega}_{t_0, \varepsilon}$ for some $\varepsilon>0$, which depends only on the geometry of $\partial\Omega_1$. This leads to a contradiction. Therefore, $v_1$ must be strictly convex in $\Omega_1$. Thus, we complete the proof of Theorem \ref{Theorem-Main}.
\end{proof}

\bibliographystyle{plain}

\bibliography{mybibliography}

\end{document}